\documentclass[12pt,oneside]{article}
\usepackage[T2A]{fontenc}
\usepackage[centertags]{amsmath}
\usepackage{amsthm}
\usepackage{amsmath}
\usepackage{amsfonts}
\usepackage{amssymb}
\usepackage{mathrsfs}
\usepackage[papersize={8.5in,11in}]{geometry}
\usepackage[colorlinks=true]{hyperref}
\newtheorem{thm}{Theorem}[section]
\newtheorem{cor}[thm]{Corollary}

\author{Asena Çetinkaya\thanks{Istanbul Kultur University, Istanbul, Turkey. Email: asnfigen@hotmail.com}~
and~Dmitrii Karp\thanks{Corresponding author.  Holon Institute of Technology, Holon, Israel. Email: dimkrp@gmail.com}}

\title {On digamma series convertible into hypergeometric series
\footnotetext{\textit{\textnormal{2020} AMS Mathematics Subject Classification:33C20,33B15,33C99}
\\ \textit{Key words and phrases: hypergeometric series, digamma series, summation formula}
}}
\date{}
\def\C{\mathbb{C}}
\def\Z{\mathbb{Z}}
\def\a{\mathbf{a}}
\def\d{\mathbf{d}}
\def\b{\mathbf{b}}
\def\cb{\mathbf{c}}
\def\cab{\mathbf{w}}
\def\eps{\varepsilon}

\def\ggamma{\boldsymbol{\gamma}}
\def\aalpha{\boldsymbol{\alpha}}
\def\bbeta{\boldsymbol{\beta}}

\def\ddelta{\boldsymbol{\delta}}

\begin{document}
\maketitle
\begin{abstract}
Series containing the digamma function arise when calculating the parametric derivatives of the hypergeometric functions and play a role in evaluation of Feynman diagrams. As these series are typically non-hypergeometric, a few instances when they are summable in terms of hypergeometric functions are  of importance.  In this paper, we convert multi-term identities for the generalized hypergeometric functions evaluated at unity into identities connecting them to the digamma series via the appropriate limiting process. The resulting formulas can be viewed as hypergeometric expressions for the $1$-norm of the gradient of the generalized hypergeometric function with respect to all its parameters and seem to have no direct  analogues in the literature.
\end{abstract}

\section{Introduction and notation}
The series containing the logarithmic derivative $\psi(z)=\Gamma'(z)/\Gamma(z)$ of Euler's gamma function $\Gamma(z)$, known as digamma or psi function,  appear in a number of contexts. First of all they may represent the parameter derivatives of hypergeometric functions,  which play an important role in several areas of mathematical physics, most notably in evaluating Feynman diagrams, see \cite{KalmKniehl,KBKMWY} and in problems involving fractional differential equations \cite{ApelMainardi}. A number of other problems from mathematical physics, where hypergeometric functions and their parameters derivatives are relevant can be found in the introductory sections of \cite{AncarGas} and \cite{Miller}.  It is thus not surprising that  the papers mentioned above as well as \cite{Cvijovic} are all published in physics or mathematical physics journals. Probably the first comprehensive collection of psi series evaluations is given in Hansen's book \cite[section 55, 360–6]{Hansen} (see also formulas numbered by $55.*$ in the addendum \cite{Borwein} by Borwein). This topic has been rather active over last two decades with a number of various results found in \cite{AncarGas,ApelMainardi,Brychkov,BrychkovGeddes,Cvijovic,Fejzullahu,Gonz-San,KangAn,Miller,RS,SB} and numerous references in these papers.  Results up to 2008 are summarized in Brychkov's book \cite{Brychkov}, see, in particular, sections 1.30.2 and 6.2. Most formulas derived in the literature convert psi series obtained by differentiating hypergeometric functions with respect to their parameters into multiple hypergeometric series, frequently of Kamp\'{e} de F\'{e}riet  type \cite{AncarGas,BrychkovGeddes,Cvijovic,Fejzullahu,KangAn,Miller,SB}.   In more special cases, however,  such series can be summed by, or rather transformed into, single hypergeometric series \cite[Chapter~3]{ApelblatBook}.  Such special cases are of importance, as the properties of the hypergeometric functions and their computation algorithms are well-developed which facilitates handling the corresponding psi series. An elegant set of formulas of this type has been established by Cvijovi\'{c} in \cite{Cvijovic} as a follow-up of the previous work by Miller \cite{Miller}, and later used by Cvijovi\'{c} and Miller to establish new reduction formulas for Kamp\'{e} de F\'{e}riet function \cite{CM}. The series handled in these references are power series whose coefficients contain the digamma function and the derivations are largely based on the Euler-Pfaff transformation for the Gauss hypergeometric function ${}_2F_1$ and the Kummer transformation for the confluent hypergeometric function  ${}_1F_1$.  Other works use series rearrangements and the hypergeometric differential equation. More recent results in this direction can be found in \cite{Gonz-San,Gonz-San23}. Our method of proof differs from all those techniques and is based on degeneration of multi-term identities for the generalized hypergeometric functions evaluated at plus or minus unity including those recently established in \cite{CKP}.  We will use this approach to transform certain rather general series containing digamma functions into combinations of univariate hypergeometric functions evaluated at the same point.  Due to significantly richer collection of transformations of hypergeometric functions at a fixed argument, our formulas have arbitrary number of free parameters unlike those due to Miller and Cvijovi\'{c} \cite{CM,Cvijovic,Miller} and those due to Gonz\'{a}lez-Santander and S\'{a}nchez Lasheras \cite{Gonz-San,Gonz-San23} as well as those  listed in \cite[6.2.1.62-81]{Brychkov}. 

To illustrate the flavour of the results anticipated in this paper let us  introduce the following  shorthand notation: 
\begin{equation}\label{eq:notation}
\begin{split}
\Gamma(\a)&=\prod\nolimits_{i=1}^{p}\Gamma(a_i),~~
(\a)_{k}=\prod\nolimits_{i=1}^{p}(a_i)_{k},~~\sin(\a)=\prod\nolimits_{i=1}^{p}\sin(a_i),
\\
\a+\beta&=(a_1+\beta,\ldots,a_p+\beta),~~\a_{[k]}=(a_1,\ldots,a_{k-1},a_{k+1},\ldots,a_p)
\end{split}
\end{equation}
for any vector $\a=(a_1,\ldots,a_p)\in\C^{p}$ (the set of complex  $p$-tuples) and  a scalar $\beta$. 
Among other things we will establish the following identities: 
\begin{multline}\label{eq:FinalMax}
\sum\limits_{k=0}^{\infty}\frac{(\a)_{k}(\pm1)^{k}}{\pi(\b)_{k}
k!}
\bigg[\sum_{b\in\b}\psi(b+k)+\psi(1+k)-\sum_{a\in\a}\psi(a+k)\bigg]
\\
\!=\!\left[\sum_{a\in\a}\cot(\pi{a})\!-\!\sum_{b\in\b}\cot(\pi{b})\right]
{}_{p}F_{p-1}\!\left(\!\left.\!\begin{array}{l}\a\\\b\end{array}\!\right|\!\pm1\!\right)
\mp\frac{(1-\b)_{1}}{\pi(1-\a)_{1}}
{}_{p+1}F_{p}\!\left(\!\left.\!\begin{array}{l}1,1,2-\b\\
2-\a\end{array}\!\right|\pm1\!\right)
\\
+\frac{\Gamma(1\!-\!\a)}{\Gamma(1\!-\!\b)}\sum\limits_{k=1}^{p-1}
\frac{[\cos({\pi}b_k)]^{\genfrac\{\}{0pt}{}{1}{0}}\Gamma(b_k\!-\!1)\Gamma(b_{k}\!-\!\b_{[k]})}{\sin({\pi}b_k)\Gamma(b_k-\a)}
{}_{p}F_{p-1}\!\left(\!\left.\!\begin{array}{l}1-b_k+\a\\2-b_k,
1-b_{k}+\b_{[k]}\end{array}\!\right|\pm1\!\right).
\end{multline}
and 
\begin{multline}\label{eq:FinalMin}
\sum\limits_{k=0}^{\infty}\frac{(\a)_{k}}{(\b)_{k}
k!}
\bigg[\sum_{b\in\b}\psi(b+k)+\psi(1+k)
-\sum_{a\in\a}\psi(a+k)\bigg]
\\
\!=\!\frac{\Gamma(\b)}{\Gamma(\a)}
\sum\limits_{k=1}^{p}\frac{\pi\Gamma(a_k)\Gamma(\a_{[k]}\!-\!a_{k})}{\sin(\pi a_k)\Gamma(\b\!-\!a_k)}
{}_{p}F_{p-1}\!\!\left(\!\!\begin{array}{l}a_{k},1-\b+a_{k}
\\
1-\a_{[k]}+a_{k}\end{array}\!\!\!\right)
-\frac{(1\!-\!\b)_{1}}{(1\!-\!\a)_{1}}
{}_{p+1}F_{p}\!\left(\!\!\begin{array}{l}1,1,2-\b\\
2-\a\end{array}\!\!\right),
\end{multline}
where here and  in the sequel ${}_{p}F_{q}(\a;\b;z)$ stands for the generalized hypergeometric function \cite[(2.1.2)]{AAR} and the omitted argument $z$ designates $z=1$.  Note that via elementary relations 
$$
\frac{d}{d{a}}(a)_n=(a)_n(\psi(a+n)-\psi(a)),~~~~\frac{d}{d{b}}\left[\frac{1}{(b)_n}\right]=\frac{1}{(b)_n}(\psi(b)-\psi(b+n))
$$
the expressions on the left hand sides of \eqref{eq:FinalMax} and \eqref{eq:FinalMin} take the form 
\begin{multline*}
\sum\limits_{k=0}^{\infty}\frac{(\a)_{k}(\pm1)^{k}}{(\b)_{k}k!}
\bigg[\sum_{b\in\b}\psi(b+k)+\psi(1+k)-\sum_{a\in\a}\psi(a+k)\bigg]
\\
=-\mathbf{1}\!\cdot\!\nabla_{\a,\b,\beta}\:{}_{p+1}F_{p}\!\left(\!\left.\!\begin{array}{l}\a,1\\\b,\beta\end{array}\!\right|\pm1\!\right)\biggl|_{ \beta=1}
+\bigg\{\sum_{b\in\b\cup\{1\}}\!\!\psi(b)-\sum_{a\in\a}\psi(a)\bigg\}
{}_{p}F_{p-1}\!\left(\!\left.\!\begin{array}{l}\a\\\b\end{array}\!\right|\pm1\!\right),
\end{multline*}
where $\mathbf{1}\!\cdot\!\nabla_{\a,\b,\beta} f(\a,\b,\beta)$ denotes the sum of the first partial derivatives of $f$ with respect to the components of $\a$, $\b$ and the variable $\beta$.  Hence, formulas  \eqref{eq:FinalMax}, \eqref{eq:FinalMin} and many further identities in this paper can be viewed as expressions for $1$-norm of the gradient of the generalized hypergeometric function with respect to all its parameters.

Another way of writing  the left hand sides of \eqref{eq:FinalMax} and \eqref{eq:FinalMin} is via Kamp\'{e} de F\'{e}riet function \cite[(28)]{SKbook}. Indeed, using 
$$
\psi(b+k)-\psi(b)=\frac{1}{b}\sum_{j=0}^{k-1}\frac{(b)_j}{(b+1)_j}
$$
and rearranging summations, we will have
\begin{multline*}
\sum\limits_{k=0}^{\infty}\frac{(\a)_{k}(\pm1)^{k}}{(\b)_{k}k!}
\bigg[\sum_{b\in\b}\psi(b+k)+\psi(1+k)-\sum_{a\in\a}\psi(a+k)\bigg]
\\
=\pm\frac{(\a)_1}{(\b)_1}\sum_{b\in\b\cup\{1\}}\frac{1}{b}
F^{p:2:1}_{p:1:0}\!\left(\left.\!\!\!\begin{array}{l}\a+1~~~:b,1~~\,:1\\\b+1,2:b+1:-\end{array}\!\right|\!\pm1,\pm1\!\right)
\\
\mp\frac{(\a)_1}{(\b)_1}\sum_{a\in\a}\frac{1}{a}F^{p:2:1}_{p:1:0}\!\left(\left.\!\!\!\begin{array}{l}\a+1~~~:a,1~~\,:1\\\b+1,2:a+1:-\end{array}\!\right|\!\pm1,\pm1\!\right)
\\
+\bigg\{\sum_{b\in\b\cup\{1\}}\!\!\psi(b)-\sum_{a\in\a}\psi(a)\bigg\}
{}_{p}F_{p-1}\!\left(\!\left.\!\begin{array}{l}\a\\\b\end{array}\!\right|\pm1\!\right).
\end{multline*}

This paper is organized as follows. In the three sections that follow we will apply degeneration process to the three identities \cite[(5.2)]{CKP}, \cite[(5.1)]{CKP} and \cite[(5.4)]{CKP}, respectively.  This results in three transformation formulas for series containing digamma functions whose combinations and particular cases yield, in particular, the above identities \eqref{eq:FinalMax} and \eqref{eq:FinalMin}. Another important particular case is given at the end of Section~4 in formula \eqref{eq:alternatingFinalMax}. Comparing the right hand sides of \eqref{eq:FinalMax} and \eqref{eq:FinalMin} as well as of other similar formulas we further get two presumably new multi-term identities for the generalized hypergeometric function evaluated at unity which are presented in formulas \eqref{eq:Fmultiterm} and \eqref{eq:FinalNoPsi}.

\section{The first transformation formula}

The formulas presented in the introduction will be derived as corollaries of several transformation formulas 
for series containing digamma functions consistent with the rising factorials. One of the denominator parameters in these formulas has to be a positive integer.  In this section, this integer represents the limiting value of the appropriate parameter difference in the identity \cite[(5.2)]{CKP} which will serve as the starting point for the degeneration process. It will be convenient to use the notation
\begin{equation}\label{eq:phi-defined}
{}_{p}\phi_{q}\!\left(\!\!\begin{array}{l}\a\\\b\!\end{array}\!\!\right)=
\frac{\Gamma(\a)}{\Gamma(\b)}{}_{p}F_{q}\!\left(\!\!\begin{array}{l}\a\\\b\!\!\end{array}\!\!\right)=\sum_{k=0}^{\infty}\frac{\Gamma(\a+k)}{\Gamma(\b+k)k!}.
\end{equation}
The main result of this section is 
\begin{thm}\label{th:first-m}
Suppose $\cb\in(\C\!\setminus\!\Z)^{p+1}$, $\d\in(\C\!\setminus\!\Z)^{p-1}$, $m\in\Z_{\ge0}$ and $\Re(\sum_{j=1}^{p-1}{d_j}-\sum_{j=1}^{p+1}(c_{j}))>-m-1$. Then
\begin{multline}\label{eq:first-m}
	\sum\limits_{k=0}^{\infty}\frac{(\cb)_{k}}{(\d)_{k}(m+k)!k!}\Big\{\sum_{d\in\d}\psi(d+k)+\psi(k+1)+\psi(1+m+k)-\sum_{c\in\cb}\psi(c+k)\Big\}
\\
=\frac{\pi}{m!}\left\{\sum_{c\in\cb}\!\cot(\pi{c})
	-\sum_{d\in\d}\!\cot(\pi{d})\right\}{}_{p+1}F_{p}\!\left(\!\!\begin{array}{l}\cb\\1+m,\d\end{array}\!\!\right)
\\
+\frac{(1-\d)_{1}}{(1-\cb)_{1}(m-1)!}
{}_{p+2}F_{p+1}\!\left(\!\!\begin{array}{l}1-m,2-\d,1,1\\2-\cb\end{array}\!\!\right)-\frac{\Gamma(1-\cb)}{\Gamma(1-\d)}
\\
\times\sum\limits_{k=1}^{p-1}\!\frac{(-1)^m\Gamma(d_{k}-m)\Gamma(d_{k})\Gamma(d_{k}-\d_{[k]})}{(d_{k}-m-1)(d_{k}-1)\Gamma(d_{k}-\cb)}
{}_{p+1}F_{p}\!\left(\!\!\begin{array}{l}1+\cb-d_{k}\\2\!-\!d_{k}, 2\!-\!d_{k}\!+\!m, 1\!+\!\d_{[k]}\!-\!d_{k}\!\!\end{array}\!\right).
\end{multline}
\end{thm}

\noindent\textbf{Proof.} We start with the  identity \cite[(5.2)]{CKP} written as follows
\begin{multline}\label{eq:expansion5.2}
	\frac{\Gamma(\b_{[1]}-b_{1})}{\Gamma(\cab-b_{1})}{}_{p+1}F_{p}\!\left(\!\!\begin{array}{l}1-\cab+b_1\\1-\b_{[1]}+b_1\!\!\end{array}\right)+\frac{\Gamma(\b_{[2]}-b_{2})}{\Gamma(\cab-b_{2})}{}_{p+1}F_{p}\!\left(\!\!\begin{array}{l}1-\cab+b_2\\1-\b_{[2]}+b_2\!\!\end{array}\right)
	\\
	=-\sum\limits_{k=3}^{p+1}\frac{\Gamma(\b_{[k]}-b_{k})}{\Gamma(\cab-b_{k})}{}_{p+1}F_{p}\!\left(\!\!\begin{array}{l}1-\cab+b_k\\1-\b_{[k]}+b_k\!\!\end{array}\right)
\end{multline}
and valid for $\b,\cab\in\C^{p+1}$ with the components of $\b$ different modulo integers. Setting $b_2=b_1+m+\eps$, $m\in\Z_{\ge0}$, using the reflection formula $\Gamma(z)\Gamma(1-z)=\pi/\sin(\pi{z})$, and the elementary relation $\sin(\pi(z+m))=(-1)^{m}\sin(\pi{z})$, for the left hand side of \eqref{eq:expansion5.2} we obtain  
\begin{multline*}
	\text{LHS of \eqref{eq:expansion5.2}}=\frac{(-1)^m}{\pi\sin(\pi\eps)}\Big(f_1(\eps)-f_2(\eps)\Big)
\\
=\frac{(-1)^m}{\sin(\pi\eps)}\Big(f_1(0)+\eps f_1'(0) -f_2(0)-\eps f_2'(0)+O(\eps^2)\Big)
\\
	=\frac{(-1)^m\eps}{\sin(\pi\eps)}\Big(f_1'(0)-f_2'(0)+O(\eps)\Big)\to\frac{(-1)^m}{\pi}\Big(f_1'(0)-f_2'(0)\Big)~\text{as}~\eps\to0,
\end{multline*}
where we employed Taylor's theorem for the functions
$$
f_1(\eps)=\frac{\sin(\pi(\cab-b_1))}{\pi\sin(\pi(\b_{[1,2]}-b_1))}{}_{p+1}\phi_{p}\!\left(\!\!\begin{array}{l}1-\cab+b_1\\1-m-\eps,1-\b_{[1,2]}+b_1\end{array}\!\!\right)
$$
and
$$
f_2(\eps)=\frac{\sin(\pi(1-\cab+b_1+\eps))}{\pi\sin(\pi(1-\b_{[1,2]}+b_1+\eps))}{}_{p+1}\phi_{p}\!\left(\!\!\begin{array}{l}1-\cab+b_1+m+\eps\\1+m+\eps,1-\b_{[1,2]}+b_1+m+\eps\!\!\end{array}\right)
$$
with $\phi$  defined in \eqref{eq:phi-defined}. The penultimate equality above is due to $f_1(0)=f_2(0)$, which is rather straightforward to verify. To compute the derivatives we will need the following two differentiation rules:
\begin{equation}\label{eq:dsin}
\frac{\partial}{\partial\eps}\left[\frac{\sin(\pi(\a+\eps))}{\sin(\pi(\d+\eps))}\right]
=\pi\frac{\sin(\pi(\a+\eps))}{\sin(\pi(\d+\eps))}\left[\sum_{a\in\a}\cot(\pi(a+\eps))-
\sum_{d\in\d}\cot(\pi(d+\eps))\right]
\end{equation}
and
\begin{equation}\label{eq:dgamma}
\frac{\partial}{\partial\eps}\left[\frac{\Gamma(\a+\eps+k)}{\Gamma(\d+\eps+k)}\right]
=\frac{\Gamma(\a+\eps+k)}{\Gamma(\d+\eps+k)}
\left[\sum_{a\in\a}\psi(a+\eps+k)-\sum_{d\in\d}\psi(d+\eps+k)\right]
\end{equation}
for any fixed vectors $\a$, $\d$.  Hence, we obtain: 
\begin{multline*}
f_1'(0)=\frac{\partial}{\partial\eps}\left[
\frac{\sin(\pi(\cab-b_1))}{\pi\sin(\pi(\b_{[1,2]}-b_1))}{}_{p+1}\phi_{p}\left(\!\!\begin{array}{l}1-\cab+b_1\\1-m-\eps,1-\b_{[1,2]}+b_1\end{array}\!\!\right)\right]_{\eps=0}
\\
=\frac{\sin(\pi(\cab-b_1))}{\pi\sin(\pi(\b_{[1,2]}-b_1))}\sum\limits_{k=0}^{m-1}
\frac{\psi(1-m+k)\Gamma(1-\cab+b_1+k)}{\Gamma(1-m+k)\Gamma(1-\b_{[1,2]}+b_1+k)k!}
\\
+\frac{\sin(\pi(\cab-b_1))}{\pi\sin(\pi(\b_{[1,2]}-b_1))}\sum\limits_{k=m}^{\infty}
\frac{\psi(1-m+k)\Gamma(1-\cab+b_1+k)}{\Gamma(1-m+k)\Gamma(1-\b_{[1,2]}+b_1+k)k!}.
\end{multline*}
In view of the reflection formula, we have 
\begin{multline*}
\frac{d}{dz}\!\left[\frac{1}{\Gamma(z)}\right]\!\!=\!\frac{d}{dz}\!\left[\frac{\sin(\pi{z})\Gamma(1-z)}{\pi}\right]
\\
\!=\!\cos(\pi{z})\Gamma(1-z)-\frac{1}{\pi}\sin(\pi{z})\psi(1-z)\Gamma(1-z)\!=\!-\frac{\psi(z)}{\Gamma(z)},
\end{multline*}
so that for $n=0,1,2,\ldots$
\begin{equation}\label{eq:psi gamma}
\frac{\psi(-n)}{\Gamma(-n)}=-\cos(\pi(-n))\Gamma(1+n)=(-1)^{n+1}n!,
\end{equation}
and, according to \eqref{eq:dsin}, the derivative $f_1'(0)$ takes the form
$$
f_1'(0)\!=\!V\!+\!\frac{\sin(\pi(\cab-b_1))}{\pi\sin(\pi(\b_{[1,2]}-b_1))}\sum\limits_{k=0}^{\infty}
\frac{\psi(1+k)\Gamma(1-\cab+b_1+m+k)}{\Gamma(1+m+k)\Gamma(1-\b_{[1,2]}+b_1+m+k)k!},
$$
where
$$
V=\frac{\sin(\pi(\cab-b_1))}{\pi\sin(\pi(\b_{[1,2]}-b_1))}\sum\limits_{k=0}^{m-1}
\frac{(-1)^{m-k}(m-k-1)!\Gamma(1-\cab+b_1+k)}{\Gamma(1-\b_{[1,2]}+b_1+k)k!}.
$$
In a similar fashion according to  \eqref{eq:dsin} and \eqref{eq:dgamma} we have:
\begin{multline*}
	f_2'(0)
 \\
 =\!\frac{\partial}{\partial\eps}\!\left[\frac{\sin(\pi(1-\cab+b_1+\eps))}{\pi\sin(\pi(1-\b_{[1,2]}+b_1+\eps))}{}_{p+1}\phi_{p}\!\left(\!\!\begin{array}{l}1-\cab+b_1+m+\eps\\1+m+\eps,1-\b_{[1,2]}+b_1+m+\eps\!\!\end{array}\right)\right]_{\eps=0}
	\\
	\!=U+\frac{\sin(\pi(\cab-b_1))}{\pi\sin(\pi(\b_{[1,2]}-b_1))}
	\sum\limits_{k=0}^{\infty}\frac{\Gamma(1-\cab+b_1+m+k)}{\Gamma(1+m+k)\Gamma(1-\b_{[1,2]}+b_1+m+k)k!}
	\\
	\times\!\Big\{\!\sum_{w\in\cab}\psi(1-w+b_1+m+k)-\!\!\sum_{b\in\b_{[1,2]}}\psi(1-b+b_1+m+k)-\psi(1+m+k)\Big\},
\end{multline*}
where
\begin{multline*}
U=\!\frac{\sin(\pi(\cab-b_1))}{\sin(\pi(\b_{[1,2]}-b_1))}\Bigg\{\sum_{b\in\b_{[1,2]}}\cot(\pi(b-b_1))
	-\sum_{w\in\cab}\cot(\pi(w-b_1))\Bigg\}
	\\
	\times{}_{p+1}\phi_{p}\left(\!\!\begin{array}{l}1-\cab+b_1+m\\1+m,1-\b_{[1,2]}+b_1+m\end{array}\!\!\right)
\end{multline*}
and we applied the elementary relations
$\sin(\pi(x-m))=(-1)^m\sin(\pi x)$, $\cot(\pi(x\pm{m}))=\cot(\pi x)$, $\cot(\pi(1-x))=-\cot(\pi{x})$.
Then, we arrive at
\begin{multline*}
	f_1'(0)-f_2'(0)\!=\!V-U
	+\frac{\sin(\pi(\cab-b_1))}{\pi\sin(\pi(\b_{[1,2]}-b_1))}\!
	\sum\limits_{k=0}^{\infty}\!\frac{\Gamma(1-\cab+b_1+m+k)}{(m+k)!\Gamma(1-\b_{[1,2]}+b_1+m+k)k!}
	\\
	\times\Big\{\psi(k+1)+\psi(1+m+k)+\sum_{b\in\b_{[1,2]}}\psi(1-b+b_1+m+k)-\sum_{w\in\cab}\psi(1-w+b_1+m+k)\Big\}.
\end{multline*}
Substituting this back into \eqref{eq:expansion5.2} for $b_2=b_1+m$ we get
\begin{multline*}
\frac{\sin(\pi(\cab-b_1))}{\pi\sin(\pi(\b_{[1,2]}-b_1))}
	\sum\limits_{k=0}^{\infty}\frac{\Gamma(1-\cab+b_1+m+k)}{\Gamma(1-\b_{[1,2]}+b_1+m+k)(m+k)!k!}
	\\
	\times\Big\{\psi(k+1)+\psi(1+m+k)+\sum\psi(1-\b_{[1,2]}+b_1+m+k)-\sum\psi(1-\cab+b_1+m+k)\Big\}
\\
=U-V-(-1)^m\pi\sum\limits_{k=3}^{p+1}\frac{\Gamma(\b_{[k]}-b_{k})}{\Gamma(\cab-b_{k})}{}_{p+1}F_{p}\!\left(\!\!\begin{array}{l}1-\cab+b_k\\1-\b_{[k]}+b_k\!\!\end{array}\right)
\end{multline*}
or, in view of definitions of $V$ and $U$:
\begin{multline*}
	\sum\limits_{k=0}^{\infty}\frac{(1-\cab+b_1+m)_{k}}{(1-\b_{[1,2]}+b_1+m)_{k}(m+k)!k!}
	\\
	\times\Big\{\psi(k+1)+\psi(1+m+k)+\sum\psi(1-\b_{[1,2]}+b_1+m+k)-\sum\psi(1-\cab+b_1+m+k)\Big\}
\\
=\frac{\pi}{m!}\!\left[\sum\limits_{b\in\b_{[1,2]}}\!\!\cot(\pi(b-b_1))
	-\sum\limits_{w\in\cab}\!\cot(\pi(w-b_1))\right]\!{}_{p+1}F_{p}\!\left(\!\!\begin{array}{l}1-\cab+b_1+m\\1+m,1-\b_{[1,2]}+b_1+m\end{array}\!\!\right)
\\
-\frac{(1-\b_{[1,2]}+b_1)_{m}}{(1-\cab+b_1)_{m}}\sum\limits_{k=0}^{m-1}
\frac{(-1)^{m-k}(m-k-1)!(1-\cab+b_1)_{k}}{(1-\b_{[1,2]}+b_1)_{k}k!}
\\
-\frac{(-1)^m\Gamma(\cab-b_1-m)}{\Gamma(\b_{[1,2]}-b_1-m)}\sum\limits_{k=3}^{p+1}\frac{\Gamma(\b_{[k]}-b_{k})}{\Gamma(\cab-b_{k})}{}_{p+1}F_{p}\!\left(\!\!\begin{array}{l}1-\cab+b_k\\1-\b_{[k]}+b_k\!\!\end{array}\right).
\end{multline*}
Writing $1-\cab+b_1+m=\cb\in\C^{p+1}$, $1-\b_{[1,2]}+b_1+m=\d\in\C^{p-1}$ we arrive at \eqref{eq:first-m} 
on application of the relation
\begin{multline}\label{eq:m-mtrick}
\frac{(1-\d)_{m}}
{(1-\cb)_{m}}\sum\limits_{k=0}^{m-1}\frac{(m-k-1)!(\cb-m)_{k}}{(-1)^{m-k}(\d-m)_{k}k!}
=\frac{(1-\d)_{1}}{(1-\cb)_{1}}\sum\limits_{j=0}^{m-1}\frac{(-1)^{j+1}(2-\d)_{j}j!}{(2-\cb)_{j}(m-1-j)!}
\\
=-\frac{(1-\d)_{1}}{(1-\cb)_{1}(m-1)!}
{}_{p+2}F_{p+1}\left(\!\!\begin{array}{l}1-m,2-\d,1,1\\2-\cb\end{array}\!\!\right),
\end{multline}
where in the first equality we used 
\begin{equation*}
\frac{(1-\d)_{m}}{(\d-m)_{k}}=(-1)^{k(p-1)}(1-\d)_{1}(2-\d)_{m-k-1}.~~~~~~~~~~~~~~~~\square
\end{equation*}

\section{The second transformation formula}
Here we will derive another hypergeometric representation for the left hand side of \eqref{eq:first-m}.  Comparing the two we will obtain a presumably new multi-term identity for the generalized hypergeometric function and the "+" case of formula \eqref{eq:FinalMax}.
\begin{thm}
Suppose $\cb\in(\C\!\setminus\!\Z)^{p+1}$, $\d\in(\C\!\setminus\!\Z)^{p-1}$, $m\in\Z_{\ge0}$ and $\Re(\sum_{j=1}^{p-1}{d_j}-\sum_{j=1}^{p+1}(c_{j}))>-m-1$. Then
\begin{multline}\label{eq:second-m}
\sum\limits_{k=0}^{\infty}\frac{(\cb)_{k}}{(\d)_{k}(m+k)!k!}
\Big\{\sum_{d\in\d}\psi(d+k)+\psi(k+1)+\psi(1+m+k)-\sum_{c\in\cb}\psi(c+k)\Big\}
\\
=\frac{\pi}{m!}\bigg[\sum_{c\in\cb_{[1]}}\cot(\pi{c})-\sum_{d\in\d_{[1]}}\cot(\pi{d})\bigg]
{}_{p+1}F_p\!\left(\begin{matrix}\cb\\1+m,\d\end{matrix}\right)
\\
+\frac{(1-\d)_{1}}{(1-\cb)_{1}(m-1)!}
{}_{p+2}F_{p+1}\left(\!\!\begin{array}{l}1-m,2-\d,1,1\\2-\cb\end{array}\!\!\right)
\\
+\frac{\Gamma(c_{1})\Gamma(d_{1})\Gamma(1-\cb_{[1]})\Gamma(1-\d_{[1]}+c_{1})}{(1-c_{1})_{m}\Gamma(d_{1}-c_{1})\Gamma(1-\cb_{[1]}+c_{1})\Gamma(1-\d_{[1]})}
{}_{p+1}F_p\!\left(\begin{matrix}c_{1},c_{1}-m,1-\d+c_{1}\\1-\cb_{[1]}+c_{1}\end{matrix}\right)
\\
-\frac{(-1)^m\Gamma(1-\cb)\sin(\pi c_{1})}{\Gamma(1-\d)\sin(\pi d_{1})}
\sum\limits_{k=2}^{p-1}
\frac{\Gamma(d_{k}-m)\Gamma(d_{k})\Gamma(d_{k}-\d_{[k]})\sin(\pi(d_{k}-d_{1}))}
{(d_{k}-m-1)(d_{k}-1)\Gamma(d_{k}-\cb)\sin(\pi(d_{k}-c_{1}))}
\\
\times{}_{p+1}F_p\!\left(\begin{matrix}1+\cb-d_{k}\\2-d_{k},2-d_{k}+m,1+\d_{[k]}-d_{k}\end{matrix}\right).
\end{multline}
\end{thm}
\textbf{Proof.} We start with identity \cite[(5.1)]{CKP} written as follows
\begin{multline}\label{eq:expansion5.1}
	\frac{\Gamma(a_2)\Gamma(\a_{[1,2]}-a_{2})}{\Gamma(1-a_{1}+a_{2})\Gamma(\b-a_{2})}
	{}_{p+1}F_p\!\left(\begin{matrix}a_{2},1-\b+a_{2}\\1-\a_{[2]}+a_{2}\end{matrix}\right)
	\\
	+\frac{\Gamma(a_3)\Gamma(\a_{[1,3]}-a_{3})}{\Gamma(1-a_{1}+a_{3})\Gamma(\b-a_{3})}
	{}_{p+1}F_p\!\left(\begin{matrix}a_{3},1-\b+a_{3}\\1-\a_{[3]}+a_{3}\end{matrix}\right)
	\\
	\!=\!
	\frac{\Gamma(\a_{[1]})}{\Gamma(1-a_{1})\Gamma(\b)}{}_{p+1}F_p\!\left(\begin{matrix}\a\\\b\end{matrix}\right)
	-\sum\limits_{k=4}^{p+1}\frac{\Gamma(a_k)\Gamma(\a_{[1,k]}-a_{k})}{\Gamma(1-a_{1}+a_{k})\Gamma(\b-a_{k})}
	{}_{p+1}F_p\!\left(\begin{matrix}a_{k},1-\b+a_{k}\\1-\a_{[k]}+a_{k}\end{matrix}\right).
\end{multline}
If we set $a_3=a_2+\eps+m$ here, then after a simple calculation similar to that in the proof of Theorem~\ref{th:first-m} we obtain
\begin{multline*}
\!\!\!\text{LHS of \eqref{eq:expansion5.1}}
\!=\!\frac{(-1)^m}{\sin(\pi\eps)}\big\{f_1(\eps)+f_2(\eps)\big\}
\!=\!\frac{(-1)^m}{\sin(\pi\eps)}\Big[f_1(0)+\eps f_1'(0) +f_2(0)+\eps f_2'(0)+O(\eps^2)\Big]
\\
=\frac{(-1)^m\eps}{\sin(\pi\eps)}\Big(f_1'(0)+f_2'(0)+O(\eps)\Big)\to\frac{(-1)^m}{\pi}\Big(f_1'(0)+f_2'(0)\Big)~\text{as}~\eps\to0,
\end{multline*}
where, in view of \eqref{eq:phi-defined},
$$
f_1(\eps)=\frac{\sin(\pi(a_1-a_2))\sin(\pi(\b-a_2))}{\pi\sin(\pi(\a_{[2,3]}-a_{2}))}
		{}_{p+1}\phi_p\!\left(\begin{matrix}a_{2},1-\b+a_{2}\\1-\eps-m,1-\a_{[2,3]}+a_{2}\end{matrix}\right),
$$
\begin{multline*}
f_2(\eps)\!=\!\frac{\sin(\pi(a_2-a_1+\eps))\sin(\pi(1-\b+a_{2}+\eps))}{\pi\sin(\pi(1-\a_{[2,3]}+a_{2}+\eps))}
\\
\times{}_{p+1}\phi_p\!\left(\begin{matrix}a_{2}+\eps+m,1-\b+a_{2}+\eps+m\\1+\eps+m, 1-\a_{[2,3]}+a_{2}+\eps+m\end{matrix}\right),
\end{multline*}
and we applied the easily verifiable identity $f_1(0)=-f_2(0)$. Then, after a simple calculation we get
\begin{multline*}
f_1'(0)=\frac{\sin(\pi(\b-a_2))}{\pi\sin(\pi(\a_{[1,2,3]}-a_{2}))}
\sum\limits_{k=0}^{m-1}\frac{(-1)^{m-k}(m-k-1)!\Gamma(a_2+k)\Gamma(1-\b+a_{2}+k)}{\Gamma(1-\a_{[2,3]}+a_{2}+k)k!}
\\
+\frac{\sin(\pi(\b-a_2))}{\pi\sin(\pi(\a_{[1,2,3]}-a_{2}))}
\sum\limits_{k=0}^{\infty}\frac{\Gamma(a_2+m+k)\Gamma(1-\b+a_{2}+m+k)\psi(1+k)}{\Gamma(1+m+k)\Gamma(1-\a_{[2,3]}+a_{2}+m+k)k!},
\end{multline*}
and, in accordance with \eqref{eq:dsin}  and \eqref{eq:dgamma},
\begin{multline*}
f_2'(0)=-\frac{\sin(\pi(\b-a_{2}))}{\sin(\pi(\a_{[1,2,3]}-a_{2}))}
\Big\{\sum_{a\in\a_{[1,2,3]}}\cot(\pi(a-a_2))-\sum_{b\in\b}\cot(\pi(b-a_2))\Big\}
\\
\times{}_{p+1}\phi_p\!\left(\begin{matrix}a_{2}+m,1-\b+a_{2}+m\\1+m,1-\a_{[2,3]}+a_{2}+m\end{matrix}\right)
\\
-\frac{\sin(\pi(\b-a_{2}))}{\pi\sin(\pi(\a_{[1,2,3]}-a_{2}))}
\sum\limits_{k=0}^{\infty}\frac{\Gamma(a_2+m+k)\Gamma(1-\b+a_{2}+m+k)}{\Gamma(1+m+k)\Gamma(1-\a_{[2,3]}+a_{2}+m+k)k!}\times
\\
\Big[\psi(a_2+m+k)+\sum_{b\in\b}\psi(1-b+a_{2}+m+k)-
\sum_{a\in\a_{[2,3]}}\psi(1-a+a_2+m+k)-\psi(1+m+k)\Big].
\end{multline*}
Substituting $(f_1'(0)+f_2'(0))/\pi$ in place of the left hand side of \eqref{eq:expansion5.1} after some tedious but elementary transformations we obtain 
\begin{multline*}
\sum\limits_{k=0}^{\infty}\frac{(a_2+m)_{k}(1-\b+a_{2}+m)_{k}}{(1+m)_{k}(1-\a_{[2,3]}+a_{2}+m)_{k}k!}
\bigg\{\sum_{a\in\a_{[2,3]}}\psi(1-a+a_2+m+k)
\\
+\psi(1+m+k)+\psi(1+k)-\psi(a_2+m+k)-\sum_{b\in\b}\psi(1-b+a_{2}+m+k)\bigg\}
\\
=
\hat{U}-\hat{V}
+\frac{(-1)^m\Gamma(\b-a_{2}-m)\Gamma(1-a_{1}+a_{2}+m)m!\Gamma(\a_{[1,2]})}{\Gamma(\a_{[1,2,3]}-a_{2}-m)(a_{2})_{m}\Gamma(1-a_{1})\Gamma(\b)}
{}_{p+1}F_p\!\left(\begin{matrix}\a\\\b\end{matrix}\right)
\\
-\frac{(-1)^m\Gamma(\b-a_{2}-m)\Gamma(1-a_{1}+a_{2}+m)}{\Gamma(\a_{[1,2,3]}-a_{2}-m)\Gamma(a_{2}+m)}
\sum\limits_{k=4}^{p+1}\frac{m!\Gamma(a_k)\Gamma(\a_{[1,k]}-a_{k})}{\Gamma(1-a_{1}+a_{k})\Gamma(\b-a_{k})}
\\
\times{}_{p+1}F_p\!\left(\begin{matrix}a_{k},1-\b+a_{k}\\1-\a_{[k]}+a_{k}\end{matrix}\right),
\end{multline*}
where
$$
\hat{U}=\pi\!\!\left[\sum\limits_{a\in\a_{[1,2,3]}}\!\!\!\cot(\pi(a-a_2))-\sum\limits_{b\in\b}\!\cot(\pi(b-a_2))\right]
\!{}_{p+1}F_p\!\left(\begin{matrix}a_{2}+m,1-\b+a_{2}+m\\1+m,1-\a_{[2,3]}+a_{2}+m\end{matrix}\right)
$$
and
$$
\hat{V}
=
\frac{(1-\a_{[2,3]}+a_2)_{m}m!}
{(1-\b+a_2)_{m}(a_{2})_{m}}
\sum\limits_{k=0}^{m-1}\frac{(-1)^{m-k}(m-k-1)!(a_2)_{k}(1-\b+a_{2})_{k}}{(1-\a_{[2,3]}+a_{2})_{k}k!}.
$$
Writing $\cb=(a_2+m,1-\b+a_2+m)\in\C^{p+1}$, $\d=(1-\a_{[2,3]}+a_2+m)\in\C^{p-1}$  and applying \eqref{eq:m-mtrick} we arrive at \eqref{eq:second-m}.
$\hfill\square$
\smallskip

Next, we derive the "+" case of \eqref{eq:FinalMax} from the above theorem.
\begin{cor}
Suppose $\cb\in(\C\!\setminus\!\Z)^{p}$, $\d\in(\C\!\setminus\!\Z)^{p-1}$ and $\Re(\sum_{j=1}^{p-1}{d_j}-\sum_{j=1}^{p}(c_{j}))>0$. Then
\begin{multline}\label{eq:final}
\sum\limits_{k=0}^{\infty}\frac{(\cb)_{k}}{(\d)_{k}k!}
\bigg\{\sum_{d\in\d}\psi(\d+k)+\psi(1+k)-\sum_{c\in\cb}\psi(\cb+k)
\bigg\}
\\
=\!\pi\!\left[\sum_{c\in\cb}\cot(\pi{c})-\sum_{d\in\d}\cot(\pi{d})\right]
{}_{p}F_{p-1}\!\left(\begin{matrix}\cb\\\d\end{matrix}\right)
\!-\!\frac{(1-\d)_{1}}{(1-\cb)_{1}}
{}_{p+1}F_p\!\left(\!\begin{matrix}1,1,2-\d\\2-\cb\end{matrix}\right)
\\
-\frac{\pi\Gamma(1-\cb)}{\Gamma(1-\d)}
\sum\limits_{k=1}^{p-1}\cot(\pi d_k)
\frac{\Gamma(d_{k}-1)\Gamma(d_{k}-\d_{[k]})}
{\Gamma(d_{k}-\cb)}
{}_{p}F_{p-1}\!\left(\begin{matrix}1-d_{k}+\cb\\2-d_{k},1-d_{k}+\d_{[k]}\end{matrix}\right).
\end{multline}
\end{cor}
\textbf{Proof}. Setting $c_1=1+m$ in \eqref{eq:second-m} we obtain
\begin{multline}\label{eq:general}
\sum\limits_{k=0}^{\infty}\frac{(\cb_{[1]})_{k}}{(\d)_{k}k!}
\left[\sum\limits_{d\in\d}\psi(d+k)+\psi(1+k)-\sum\limits_{c\in\cb_{[1]}}\psi(c+k)
\right]
\\
=\pi\left[\sum\limits_{c\in\cb_{[1]}}\cot(\pi{c})-\sum\limits_{d\in\d_{[1]}}\cot(\pi{d})\right]
{}_{p}F_{p-1}\!\left(\begin{matrix}\cb_{[1]}\\\d\end{matrix}\right)-m!\hat{V}
\\
+\frac{(-1)^mm!\Gamma(d_{1})\Gamma(1-\cb_{[1]})\Gamma(1-\d_{[1]}+m+1)}{\Gamma(d_{1}-m-1)\Gamma(1-\cb_{[1]}+m+1)\Gamma(1-\d_{[1]})}
{}_{p+1}F_p\!\left(\begin{matrix}m+1,1,1-\d+c_{1}\\1-\cb_{[1]}+c_{1}\end{matrix}\right)
\\
+\frac{\Gamma(d_{1})\Gamma(1-\cb_{[1]})}{\Gamma(1-\d_{[1]})}
\sum\limits_{k=2}^{p-1}
\frac{\Gamma(d_{k}-1)\Gamma(d_{k})\Gamma(1-d_{k})\Gamma(d_{k}-\d_{[1,k]})}
{\Gamma(d_{k}-\cb_{[1]})\Gamma(1-d_{k}+d_{1})}
\\
\times\!{}_{p}F_{p-1}\!\left(\begin{matrix}1-d_{k}+\cb_{[1]}\\2-d_{k},1-d_{k}+\d_{[k]}\end{matrix}\right),
\end{multline}
where 
$$
\hat{V}=\frac{(1-\d)_{m}}
{(1-\cb)_{m}}\sum\limits_{k=0}^{m-1}\frac{(-1)^{m-k}(m-k-1)!(\cb-m)_{k}}{(\d-m)_{k}k!}
=\frac{(1-\d)_{1}}{m!(1-\cb{[1]})_{1}}\sum\limits_{j=0}^{m-1}\frac{j!(2-\d)_{j}}{(2-\cb_{[1]})_{j}}
$$
(we used that  $c_1=m+1$ in the last equality). It is easy to see that for $c_1=m+1$ in view of  $(a+m)_{j}=(a)_{m+j}/(a)_{m}$ we will have
\begin{multline}\label{eq:m-aux}
-m!\hat{V}+\frac{(-1)^mm!\Gamma(d_{1})\Gamma(1-\cb_{[1]})\Gamma(1-\d_{[1]}+m+1)}{\Gamma(d_{1}-m-1)\Gamma(1-\cb_{[1]}+m+1)\Gamma(1-\d_{[1]})}
\\
\times{}_{p+1}F_p\!\left(\begin{matrix}m+1,1,1-\d+c_{1}\\1-\cb_{[1]}+c_{1}\end{matrix}\right)
=-\frac{(1-\d)_{1}}{(1-\cb_{[1]})_{1}}
{}_{p+1}F_p\!\left(\begin{matrix}1,1,2-\d\\2-\cb_{[1]}\end{matrix}\right).
\end{multline}
Substituting \eqref{eq:m-aux} into \eqref{eq:second-m} and changing  $\cb_{[1]}\to\cb$ we  get
\begin{multline*} 
\sum\limits_{k=0}^{\infty}\frac{(\cb)_{k}}{(\d)_{k}k!}
\left[\sum\limits_{d\in\d}\psi(d+k)+\psi(1+k)-\sum\limits_{c\in\cb}\psi(\cb+k)
\right]
\\
=\pi\left[\sum\limits_{c\in\cb}\cot(\pi{c})-\sum\limits_{d\in\d_{[1]}}\cot(\pi{d})\right]
{}_{p}F_{p-1}\!\left(\begin{matrix}\cb\\\d\end{matrix}\right)
-\frac{(1-\d)_{1}}{(1-\cb)_{1}}
{}_{p+1}F_p\!\left(\begin{matrix}1,1,2-\d\\2-\cb\end{matrix}\right)
\\
+\frac{\pi\Gamma(1-\cb)}{\Gamma(1-\d)}
\sum\limits_{k=2}^{p-1}
\frac{\sin[\pi(d_{k}-d_{1})]\Gamma(d_{k}-1)\Gamma(d_{k}-\d_{[k]})}
{\sin(\pi d_{1})\sin(\pi d_k)\Gamma(d_{k}-\cb)}
{}_{p}F_{p-1}\!\left(\begin{matrix}1-d_{k}+\cb\\2-d_{k},1-d_{k}+\d_{[k]}\end{matrix}\right).
\end{multline*}
Using  
$$
\frac{\sin(\pi(d_{k}-d_{1}))}{\sin(\pi d_{1})\sin(\pi d_k)}=\cot(\pi d_{1})-\cot(\pi d_k)
$$
and expressing from \cite[(5.3)]{CKP} 
\begin{multline*}
\sum\limits_{k=2}^{p-1}
\frac{\Gamma(d_{k}-1)\Gamma(d_{k}-\d_{[k]})}
{\Gamma(d_{k}-\cb)}
{}_{p}F_{p-1}\!\left(\begin{matrix}1-d_{k}+\cb\\2-d_{k},1-d_{k}+\d_{[k]}\end{matrix}\right)
\\
=-\frac{\Gamma(1-\d)}{\Gamma(1-\cb)}{}_{p}F_{p-1}\!\left(\begin{matrix}\cb\\\d\end{matrix}\right) 
-\frac{\Gamma(d_{1}-1)\Gamma(d_{1}-\d_{[1]})}
{\Gamma(d_{1}-\cb)}
{}_{p}F_{p-1}\!\left(\begin{matrix}1-d_{1}+\cb\\2-d_{1},1-d_{1}+\d_{[1]}\end{matrix}\right),
\end{multline*}
we finally arrive at \eqref{eq:final}. $\hfill\square$

Equating \eqref{eq:first-m} and \eqref{eq:second-m} after some rearrangements and trigonometric manipulations we obtain
\begin{multline*}
\frac{\Gamma(a_{1})\Gamma(d_{1})\Gamma(1-\d_{[1]}+a_{1})}{(1-a_{1})_{m}\Gamma(d_{1}-a_{1})\Gamma(1-\a_{[1]}+a_{1})}
{}_{p+1}F_p\!\left(\begin{matrix}a_{1},a_{1}-m,1-\d+a_{1}\\1-\a_{[1]}+a_{1}\end{matrix}\right)
\!+\!\frac{\Gamma(d_{1}-\d_{[1]})}{\Gamma(d_{1}-\a)}
\\
\times\!\frac{\Gamma(d_{1}-m-1)}{(-1)^m\Gamma(1-d_1)}\Gamma(1-a_1)\Gamma(d_{1}-1)
{}_{p+1}F_{p}\!\left(\!\!\begin{array}{l}1-d_{1}+\a\\2-d_{1}, 2-d_{1}+m, 1-d_{1}+\d_{[1]}\!\!\end{array}\right)
\\
-\frac{\pi\Gamma(1-\d_{[1]})}{m!\Gamma(1-\a_{[1]})}\left(\cot(\pi a_1)
	-\cot(\pi d_1)\right){}_{p+1}F_{p}\left(\!\!\begin{array}{l}\a\\1+m,\d\end{array}\!\!\right)
\\
=-(-1)^m
\sum\limits_{k=2}^{p-1}
\frac{\Gamma(d_{k}-m-1)\Gamma(d_{k}-\d_{[k]})\Gamma(d_1)\Gamma(1-a_1)\Gamma(1-d_k+a_1)}
{\Gamma(2-d_{k})\Gamma(a_{1}-d_{1})\Gamma(1-a_{1}+d_{1})\Gamma(d_{k}-\a_{[1]})}
\\
\times{}_{p+1}F_{p}\!\left(\!\!\begin{array}{l}1-d_{k}+\a\\2-d_{k}, 2-d_{k}+m, 1-d_{k}+\d_{[k]}\!\!\end{array}\right).
\end{multline*}
This leads to 
\begin{multline*}
\frac{\Gamma(a_{1})\Gamma(d_{1})\Gamma(1-\d_{[1]}+a_{1})}{(1-a_{1})_{m}\Gamma(d_{1}-a_{1})\Gamma(1-\a_{[1]}+a_{1})}
{}_{p+1}F_p\!\left(\begin{matrix}a_{1},a_{1}-m,1-\d+a_{1}\\1-\a_{[1]}+a_{1}\end{matrix}\right)
\\
-\frac{\pi\Gamma(1-\d_{[1]})}{m!\Gamma(1-\a_{[1]})}\left(\cot(\pi a_1)
	-\cot(\pi d_1)\right){}_{p+1}F_{p}\left(\!\!\begin{array}{l}\a\\1+m,\d\end{array}\!\!\right)
\\
=-(-1)^m
\sum\limits_{k=1}^{p-1}
\frac{\Gamma(d_1)\Gamma(d_{k}-m-1)\Gamma(d_{k}-\d_{[k]})\Gamma(1-a_1)\Gamma(1-d_k+a_1)}
{\Gamma(2-d_{k})\Gamma(a_{1}-d_{1})\Gamma(1-a_{1}+d_{1})\Gamma(d_{k}-\a_{[1]})}
\\
\times{}_{p+1}F_{p}\!\left(\!\!\begin{array}{l}1-d_{k}+\a\\2-d_{k}, 2-d_{k}+m, 1-d_{k}+\d_{[k]}\!\!\end{array}\right).
\end{multline*}

Rearranging the above formula, we obtain the following presumably new  $p+1$-term identity.
\begin{cor}
the following identity holds for $m\in\mathbb{N}$ and  the values of parameters making the hypergeometric series involved finite\emph{:}
\begin{multline}\label{eq:Fmultiterm}
-\frac{\Gamma(1-\d)}{m!\Gamma(1-\a)}{}_{p+1}F_{p}\left(\!\!\begin{array}{l}\a\\1+m,\d\end{array}\!\!\right)
+\frac{\sin(\pi a_{1})\Gamma(a_{1})\Gamma(1-\d+a_{1})}{\pi\Gamma(1-a_{1}+m)\Gamma(1-\a_{[1]}+a_{1})}
\\
\times{}_{p+1}F_p\!\left(\begin{matrix}a_{1},a_{1}-m,1-\d+a_{1}\\1-\a_{[1]}+a_{1}\end{matrix}\right)
=(-1)^m
\sum\limits_{k=1}^{p-1}
\frac{\Gamma(d_{k}-m-1)\Gamma(d_{k}-\d_{[k]})}{\Gamma(2-d_{k})\Gamma(d_{k}-\a)}
\\
\times\frac{\sin(\pi a_{1})}{\sin(\pi(d_k-a_1))}{}_{p+1}F_{p}\!\left(\!\!\begin{array}{l}1-d_{k}+\a\\2-d_{k}, 2-d_{k}+m, 1-d_{k}+\d_{[k]}\!\!\end{array}\right).
\end{multline}
\end{cor}
This formula can be viewed as the degenerate case of \cite[(5.3)]{CKP} when one of the bottom parameters is positive integer.  We emphasize that computing the corresponding limit in \cite[(5.3)]{CKP} will lead to psi series, not to the above formula which appears to be significantly  less trivial.  

\section{The third transformation formula}

Suppose $0\le{n}\le{p}$, $0\le{s}\le{p}$ are integers that satisfy $s+n\ge{p}$ and $\a\in\C^{n}$,  $\b\in\C^{s}$, $\cb\in\C^{p-n}$, $\d\in\C^{p-s}$ are complex vectors that satisfy
$$
\Re\Bigl(\sum\nolimits_{j=1}^{n}a_j+\sum\nolimits_{j=1}^{p-n}c_j-\sum\nolimits_{j=1}^{s}b_j-\sum\nolimits_{j=1}^{p-s}d_j\Bigr)>0.
$$
Under these conditions we have recently established the identity \cite[Theorem~5.1]{CKP}
\begin{subequations}\label{eq:Newsummation}
\begin{multline}\label{eq:CKPTh5.1}
\sum\limits_{k=1}^{s}\frac{A_{k}}{b_{k}}
{}_{p+1}F_{p}\!\left(\!\left.\!\begin{array}{l}1-\a+b_{k},1-\cb+b_{k},b_{k}\\1-\b_{[k]}+b_{k},1-\d+b_{k},b_{k}+1\end{array}\!\right|(-1)^{p-s-n}\!\right)
\\
+\sum\limits_{k=1}^{n}\!\frac{B_{k}}{(1-a_{k})}
{}_{p+1}F_{p}\!\left(\!\left.\!\begin{array}{l}1+\b-a_{k},1+\d-a_{k},1-a_{k}\\1+\a_{[k]}-a_{k},1+\cb-a_{k},2-a_{k}\end{array}\!\right|(-1)^{p-s-n}\!\right)
\!=\!\frac{\Gamma(1-\a)\Gamma(\b)}{\Gamma(\cb)\Gamma(1-\d)},
\end{multline}
where
\begin{equation}\label{eq:constants}
A_{k}=\frac{\Gamma(\b_{[k]}-b_{k})\Gamma(1-\a+b_k)}{\Gamma(\cb-b_{k})\Gamma(1-\d+b_{k})},~~
B_k=\frac{\Gamma(a_k-\a_{[k]})\Gamma(1+\b-a_{k})}{\Gamma(a_{k}-\d)\Gamma(1+\cb-a_{k})},
\end{equation}
\end{subequations}
and it is assumed that no numerator gamma factor becomes infinite.  If the latter condition is violated it is still possible to compute the corresponding limit leading naturally to psi series. Using this approach we will establish the following
\begin{thm}
Suppose $0\le{n}\le{p}$, $0\le{s}\le{p-1}$, $\kappa:=s-n\ge-2$ and 
\begin{align*}
&\a_1\in(\C\!\setminus\!\Z)^{n},~~ \a_2\in\C^{p-n},~~ \a=(\a_1,\a_2)\in\C^{p}, 
\\
&\b_1\in(\C\!\setminus\!\Z)^{s},~~ \b_2\in\C^{p-1-s},~~ \b=(\b_1,\b_2)\in\C^{p-1}. 
\end{align*}
Assume further that $\Re\big(\sum_{b\in\b}b-\sum_{a\in\a}a\big)>0$. Then
\begin{multline}\label{eq:alternatingFinal}
\sum\limits_{k=0}^{\infty}\frac{(\a)_{k}[(-1)^{\kappa}]^{k}}{\pi(\b)_{k}
k!}
\bigg[\!\sum_{b\in\b}\psi(b+k)+\psi(1+k)
-\sum_{a\in\a}\psi(a+k)\bigg]
\\
\!=\!\Bigg\{\sum_{a\in\a_1}\cot(\pi{a})\!-\!\sum_{b\in\b_1}\cot(\pi{b})\Bigg\}{}_{p}F_{p-1}\!\left(\!\left.\!\!\begin{array}{l}\a\\\b\end{array}\!\right|(-1)^{\kappa}\!\right)
\\
-\frac{(-1)^{\kappa}(1-\b)_{1}}{\pi(1-\a)_{1}}
{}_{p+1}F_{p}\!\left(\!\left.\!\!\begin{array}{l}1,1,2-\b\\
2-\a\end{array}\!\right|(-1)^{\kappa}\!\right)+\frac{\Gamma(1-\a_{1})\Gamma(\b_{2})}{\Gamma(1-\b_{1})\Gamma(\a_{2})}
\\
\times\!\sum\limits_{a_k\in\a_{2}}\!\frac{\Gamma(a_k)\Gamma(\a_{2[k]}-a_{k})\Gamma(1-\b_1+a_k)}{\sin(\pi a_k)\Gamma(\b_2-a_k)\Gamma(1-\a_1+a_{k})}
{}_{p}F_{p-1}\!\left(\!\left.\!\!\begin{array}{l}a_{k},1-\b+a_{k}
\\
1-\a_{[k]}+a_{k}\end{array}\!\!\right|(-1)^{\kappa}\!\right)
\\
+\frac{\Gamma(1-\a_{1})\Gamma(\b_{2})}{\Gamma(1-\b_{1})\Gamma(\a_{2})}\sum\limits_{b_k\in\b_1}\!\!
\frac{\Gamma(b_k-1)\Gamma(b_{k}-\b_{1[k]})\Gamma(1-b_k+\a_{2})}{\sin({\pi}b_k)\Gamma(b_k-\a_{1})\Gamma(1-b_k+\b_{2})}
\\
\times{}_{p}F_{p-1}\!\left(\!\left.\!\begin{array}{l}1-b_k+\a\\2-b_k,
1-b_{k}+\b_{[k]}\end{array}\!\right|\!(-1)^{\kappa}\!\right).
\end{multline}
\end{thm}

\textbf{Proof.} Setting  $b_2=b_1+m+\eps$ the first two terms on the left hand side of \eqref{eq:CKPTh5.1}  with coefficients $A_1$, $A_2$ take the form (put $\varkappa=p-s-n$ for brevity):
\begin{multline*}
\frac{\Gamma(\eps+m)\Gamma(\b_{[1,2]}-b_{1})\Gamma(1-\a+b_1)}{b_1\Gamma(\cb-b_1)\Gamma(1-\d+b_{1})}
\\
\times{}_{p+1}F_{p}\!\left(\!\left.\!\begin{array}{l}1-\a+b_{1},1-\cb+b_{1},b_{1}\\1-m-\eps,1-\b_{[1,2]}+b_{1},1-\d+b_{1},b_{1}+1\end{array}\!\right|(-1)^{\varkappa}\!\right)
\\
+\frac{\Gamma(-\eps-m)\Gamma(\b_{[1,2]}-b_{1}-m-\eps)\Gamma(1-\a+b_1+m+\eps)}{(b_1+m+\eps)\Gamma(\cb-b_1-m-\eps)\Gamma(1-\d+b_{1}+m+\eps)}
\\
\times\!{}_{p+1}F_{p}\!\left(\!\left.\!\!\begin{array}{l}1\!-\!\a\!+\!b_{1}\!+\!m\!+\!\eps,1\!-\!\cb\!+\!b_{1}\!+\!m\!+\!\eps,b_{1}\!+\!m\!+\!\eps\\1\!+\!m\!+\!\eps,1\!-\!\b_{[1,2]}\!+\!b_{1}\!+\!m\!+\!\eps,1\!-\!\d\!+\!b_{1}\!+\!m\!+\!\eps,b_{1}\!+\!m\!+\!\eps\!+\!1\end{array}\!\!\right|\!(-1)^{\varkappa}\!\!\right).
\end{multline*}
Simple transformations of this expression leads to  
\begin{equation*}
	\frac{(-1)^m}{\pi^{\varkappa+1}\sin(\pi\eps)}\big\{f_1(\eps)
	-f_2(\eps)\big\}\!=\!\frac{(-1)^m}{\pi^{\varkappa+1}\sin(\pi\eps)}\Big[f_1(0)+\eps f_1'(0) -f_2(0)-\eps f_2'(0)+O(\eps^2)\Big]
\end{equation*}
as $\eps\to0$, where $\phi$ is defined in \eqref{eq:phi-defined} and
$$
f_1(\eps)\!=\!\frac{\sin(\pi(\cb-b_1))}{\sin(\pi(\b_{[1,2]}-b_{1}))}
		{}_{p+1}\phi_{p}\!\left(\!\left.\!\!\begin{array}{l}1-\a+b_{1},1-\cb+b_{1},b_{1}\\1-m-\eps,1-\b_{[1,2]}+b_{1},1-\d+b_{1},b_{1}+1\end{array}\!\!\right|\!(-1)^{\varkappa}\!\right),
$$
\begin{multline*}
f_2(\eps)=\frac{(-1)^{m\varkappa}\sin(\pi(1-\cb+b_{1}+\eps))}{\sin(\pi(1-\b_{[1,2]}+b_{1}+\eps))} 
\\ 
\times{}_{p+1}\phi_{p}\!\!\left(\!\left.\!\!\begin{array}{l}1-\a+b_{1}+m+\eps,1-\cb+b_{1}+m+\eps,b_{1}+m+\eps\\1\!+\!m\!+\!\eps,1\!-\!\b_{[1,2]}\!+\!b_{1}\!+\!m\!+\!\eps,1\!-\!\d\!+\!b_{1}\!+\!m\!+\!\eps,b_{1}\!+\!m\!+\!\eps\!+\!1\end{array}\!\!\right|\!(-1)^{\varkappa}\!\right).
\end{multline*}
Note that $f_1(0)=f_2(0)$, because $\sin(\pi(1-\cb+b_{1}))=\sin(\pi(\cb-b_1))$ and 
\begin{multline*}
{}_{p+1}\phi_{p}\!\left(\!\left.\!\begin{array}{l}1-\a+b_{1},1-\cb+b_{1},b_{1}\\1-m,1-\b_{[1,2]}+b_{1},1-\d+b_{1},b_{1}+1\end{array}\!\right|(-1)^{\varkappa}\!\right)
\\
=\sum\limits_{k=0}^{\infty}\frac{\Gamma(1-\a+b_1+k)\Gamma(1-\cb+b_{1}+k)\Gamma(b_1+k)[(-1)^{\varkappa}]^k}{\Gamma(1-m+k)\Gamma(1-\b_{[1,2]}+b_{1}+k)\Gamma(1-\d+b_{1}+k)\Gamma(b_1+1+k)k!}
\\
=\sum\limits_{k=m}^{\infty}\frac{\Gamma(1-\a+b_1+k)\Gamma(1-\cb+b_{1}+k)\Gamma(b_1+k)[(-1)^{\varkappa}]^k}{\Gamma(1-m+k)\Gamma(1-\b_{[1,2]}+b_{1}+k)\Gamma(1-\d+b_{1}+k)\Gamma(b_1+1+k)k!}
\\
=\sum\limits_{k=0}^{\infty}\frac{\Gamma(1-\a+b_1+m+k)\Gamma(1-\cb+b_{1}+m+k)\Gamma(b_1+m+k)[(-1)^{\varkappa}]^k}{k!\Gamma(1\!-\!\b_{[1,2]}\!+\!b_{1}\!+\!m\!+\!k)\Gamma(1\!-\!\d\!+\!b_{1}\!+\!m\!+\!k)\Gamma(b_1\!+\!1\!+\!m\!+\!k)(m+k)!}
\\
=(-1)^{m\varkappa}{}_{p+1}\phi_{p}\!\left(\!\left.\!\begin{array}{l}1-\a+b_{1}+m,1-\cb+b_{1}+m,b_{1}+m\\1+m,1-\b_{[1,2]}+b_{1}+m,1-\d+b_{1}+m,b_{1}+m+1\end{array}\!\right|(-1)^{\varkappa}\!\right).
\end{multline*}
Hence,
\begin{multline*}
\frac{(-1)^m}{\pi^{\varkappa+1}\sin(\pi\eps)}\big\{f_1(\eps)
	-f_2(\eps)\big\}
 \\
 =\frac{(-1)^m\eps}{\pi^{\varkappa+1}\sin(\pi\eps)}\Big(f_1'(0)-f_2'(0)+O(\eps)\Big)\to\frac{(-1)^m}{\pi^{\varkappa+2}}\Big(f_1'(0)-f_2'(0)\Big)~\text{as}~\eps\to0.
\end{multline*}
In view of $\frac{d}{d\varepsilon}\left[\Gamma(1-m-\varepsilon)\right]^{-1}=\psi(1-m-\varepsilon)/\Gamma(1-m-\varepsilon)$ and formula \eqref{eq:psi gamma}, we get by separating summations over $k=0,\ldots,m-1$ and $k=m,m+1,\ldots$: 
\begin{multline*}
	f_1'(0)=\frac{\sin(\pi(\cb-b_1))}{\sin(\pi(\b_{[1,2]}-b_{1}))}
\\
\times\!\sum\limits_{k=0}^{m-1}\frac{(-1)^{m-k}(m-k-1)![(-1)^{\varkappa}]^k\Gamma(1-\a+b_1+k)\Gamma(1-\cb+b_{1}+k)\Gamma(b_1+k)}{\Gamma(1-\b_{[1,2]}+b_{1}+k)\Gamma(1-\d+b_{1}+k)\Gamma(b_1+1+k)k!}
	\\+ \frac{\sin(\pi(\cb-b_1))}{\sin(\pi(\b_{[1,2]}-b_{1}))}
 \\
\times\!\sum\limits_{k=0}^{\infty}\!\frac{[(-1)^{\varkappa}]^{k+m}\Gamma(1\!-\!\a\!+\!b_1\!+\!m\!+\!k)\Gamma(1\!-\!\cb\!+\!b_{1}\!+\!m\!+\!k)\Gamma(b_1+m+k)\psi(1+k)}
{\Gamma(1\!-\!\b_{[1,2]}!+!b_{1}\!+\!m\!+\!k)\Gamma(1\!-\!\d\!+\!b_{1}\!+\!m\!+\!k)\Gamma(b_1+1+m+k)\Gamma(1+m+k)k!}
\end{multline*}
or
\begin{multline*}
	f_1'(0)=\frac{\pi^{\varkappa+2}\Gamma(1-\a+b_1)\Gamma(\b_{[1,2]}-b_{1})}{\Gamma(\cb-b_1)\Gamma(1-\d+b_{1})}
 \\
	\times\!\sum\limits_{k=0}^{m-1}\frac{(-1)^{m-k}(m-k-1)![(-1)^{\varkappa}]^k(1-\a+b_1)_{k}(1-\cb+b_{1})_{k}}{(1-\b_{[1,2]}+b_{1})_{k}(1-\d+b_{1})_{k}(b_1+k)k!}
\\
	+\frac{\pi^{\varkappa+2}\Gamma(\b_{[1,2]}-b_{1}-m)\Gamma(1-\a+b_1+m)}{\Gamma(\cb-b_1-m)\Gamma(1-\d+b_{1}+m)}
 \\
\times\!\sum\limits_{k=0}^{\infty}\!\frac{(-1)^{\varkappa{k}}(1-\a+b_1+m)_{k}(1-\cb+b_{1}+m)_{k}\psi(1+k)}{(1-\b_{[1,2]}+b_{1}+m)_{k}(1-\d+b_{1}+m)_{k}(b_1+m+k)(m+k)!k!}.
\end{multline*}
Similarly,  by applying  \eqref{eq:dsin} and \eqref{eq:dgamma} and in view of $\sin(\pi(1-\alpha))=\sin(\pi\alpha)$ and $\cot(\pi(1-\alpha))=-\cot(\pi\alpha)$, we get
\begin{multline*}
f_2'(0)=\frac{(-1)^{m\varkappa}\pi\sin(\pi(\cb-b_{1}))}{\sin(\pi(\b_{[1,2]}-b_{1}))}
\Biggl[\sum\limits_{b\in\b_{[1,2]}}\!\!\!\cot(\pi(b-b_{1}))-\sum\limits_{c\in\cb}\!\!\cot(\pi(c-b_{1}))\Biggr]
\\
\times{}_{p+1}\phi_{p}\!\left(\!\left.\!\begin{array}{l}1-\a+b_{1}+m,1-\cb+b_{1}+m,b_{1}+m\\1+m,1-\b_{[1,2]}+b_{1}+m,1-\d+b_{1}+m,b_{1}+m+1\end{array}\!\right|(-1)^{\varkappa}\!\right)
		\\
+\frac{(-1)^{m\varkappa}\sin(\pi(\cb-b_{1}))}{\sin(\pi(\b_{[1,2]}-b_{1}))}
\\
\times\!\sum\limits_{k=0}^{\infty}\frac{[(-1)^{\varkappa}]^k\Gamma(1\!-\!\a\!+\!b_1\!+\!m\!+\!k)\Gamma(1\!-\!\cb\!+\!b_{1}\!+\!m\!+\!k)\Gamma(b_1+m+k)\Psi}{\Gamma(1\!-\!\b_{[1,2]}\!+\!b_{1}\!+\!m\!+\!k)\Gamma(1\!-\!\d\!+\!b_{1}\!+\!m\!+\!k)\Gamma(b_1+1+m+k)(m+k)!k!}
		\\
=\!\frac{\pi^{\varkappa+3}\Gamma(\b_{[1,2]}-b_{1}-m)\Gamma(1-\a+b_{1}+m)}{\Gamma(\cb-b_{1}-m)\Gamma(1-\d+b_{1}+m)(b_{1}+m)m!}
\Biggl[\sum\limits_{b\in\b_{[1,2]}}\!\!\!\cot(\pi(b-b_{1}))-\sum\limits_{c\in\cb}\!\!\cot(\pi(c-b_{1}))\Biggr]
\\
\times\!{}_{p+1}F_{p}\!\left(\!\left.\!\begin{array}{l}1-\a+b_{1}+m,1-\cb+b_{1}+m,b_{1}+m\\1+m,1-\b_{[1,2]}+b_{1}+m,1-\d+b_{1}+m,b_{1}+m+1\end{array}\!\right|(-1)^{\varkappa}\!\right)
\\
+\frac{\pi^{\varkappa+2}\Gamma(\b_{[1,2]}-b_{1}-m)\Gamma(1-\a+b_{1}+m)}{\Gamma(\cb-b_{1}-m)\Gamma(1-\d+b_{1}+m)}
\\
\times\!\sum\limits_{k=0}^{\infty}\frac{[(-1)^{\varkappa}]^k(1-\a+b_1+m)_{k}(1-\cb+b_{1}+m)_{k}\Psi}{(1-\b_{[1,2]}+b_{1}+m)_{k}(1-\d+b_{1}+m)_{k}(b_1+m+k)(m+k)!k!},
\end{multline*}
where
\begin{multline*}
\Psi=\sum\limits_{a\in\a}\psi(1-a+b_1+m+k)+\sum\limits_{c\in\cb}\psi(1-c+b_{1}+m+k)+1/(b_1+m+k)
\\
-\!\!\!\sum\limits_{b\in\b_{[1,2]}}\!\!\!\!\psi(1-b+b_{1}+m+k)-\sum\limits_{d\in\d}\!\psi(1-d+b_{1}+m+k)-\psi(1+m+k).
\end{multline*}
Combining the above expressions we obtain
\begin{multline*}
\frac{(-1)^m}{\pi^{\varkappa+2}}\Big(f_1'(0)-f_2'(0)\Big)=
\frac{(-1)^m\Gamma(1-\a+b_1)\Gamma(\b_{[1,2]}-b_{1})}{\Gamma(\cb-b_1)\Gamma(1-\d+b_{1})}
\\
\times\!\sum\limits_{k=0}^{m-1}\frac{(-1)^{m-k}(m-k-1)![(-1)^{\varkappa}]^k(1-\a+b_1)_{k}(1-\cb+b_{1})_{k}}{(1-\b_{[1,2]}+b_{1})_{k}(1-\d+b_{1})_{k}(b_1+k)k!}
\\
-\frac{\pi(-1)^m\Gamma(\b_{[1,2]}-b_{1}-m)\Gamma(1-\a+b_{1}+m)}{\Gamma(\cb-b_{1}-m)\Gamma(1-\d+b_{1}+m)(b_{1}+m)m!}
\Biggl[\sum\limits_{b\in\b_{[1,2]}}\!\!\!\cot(\pi(b-b_{1}))-\sum\limits_{c\in\cb}\!\!\cot(\pi(c-b_{1}))\Biggr]
\\
\times{}_{p+1}F_{p}\!\left(\!\left.\!\begin{array}{l}1-\a+b_{1}+m,1-\cb+b_{1}+m,b_{1}+m\\1+m,1-\b_{[1,2]}+b_{1}+m,1-\d+b_{1}+m,b_{1}+m+1\end{array}\!\right|(-1)^{\varkappa}\!\right)
\\
-\frac{(-1)^m\Gamma(\b_{[1,2]}-b_{1}-m)\Gamma(1-\a+b_{1}+m)}{\Gamma(\cb-b_{1}-m)\Gamma(1-\d+b_{1}+m)}
\\
\times\!\sum\limits_{k=0}^{\infty}\frac{[(-1)^{\varkappa}]^k(1-\a+b_1+m)_{k}(1-\cb+b_{1}+m)_{k}\Psi}{(1-\b_{[1,2]}+b_{1}+m)_{k}(1-\d+b_{1}+m)_{k}(b_1+m+k)(m+k)!k!}.
\end{multline*}
Substitution back into \eqref{eq:Newsummation}  then yields:
\begin{multline*}
\frac{\Gamma(1-\a+b_1)\Gamma(\b_{[1,2]}-b_{1})}{(-1)^m\Gamma(\cb-b_1)\Gamma(1-\d+b_{1})}
	\sum\limits_{k=0}^{m-1}\frac{(m-k-1)![(-1)^{\varkappa}]^k(1-\a+b_1)_{k}(1-\cb+b_{1})_{k}}{(-1)^{m-k}(1-\b_{[1,2]}+b_{1})_{k}(1-\d+b_{1})_{k}(b_1+k)k!}
\\
-\frac{\pi(-1)^m\Gamma(\b_{[1,2]}-b_{1}-m)\Gamma(1-\a+b_{1}+m)}{\Gamma(\cb-b_{1}-m)\Gamma(1-\d+b_{1}+m)(b_{1}+m)m!}
\Biggl[\sum\limits_{b\in\b_{[1,2]}}\!\!\!\cot(\pi(b-b_{1}))-\sum\limits_{c\in\cb}\!\!\cot(\pi(c-b_{1}))\Biggr]
\\
\times{}_{p+1}F_{p}\!\left(\!\left.\!\begin{array}{l}1-\a+b_{1}+m,1-\cb+b_{1}+m,b_{1}+m\\1+m,1-\b_{[1,2]}+b_{1}+m,1-\d+b_{1}+m,b_{1}+m+1\end{array}\!\right|(-1)^{\varkappa}\!\right)
\\
+\frac{(-1)^m\Gamma(\b_{[1,2]}-b_{1}-m)\Gamma(1-\a+b_{1}+m)}{\Gamma(\cb-b_{1}-m)\Gamma(1-\d+b_{1}+m)}
\\
\times\sum\limits_{k=0}^{\infty}\frac{[(-1)^{\varkappa}]^k(1-\a+b_1+m)_{k}(1-\cb+b_{1}+m)_{k}\Psi}{(1-\b_{[1,2]}+b_{1}+m)_{k}(1-\d+b_{1}+m)_{k}(b_1+m+k)(m+k)!k!}
\\
=\frac{\Gamma(1-\a)\Gamma(\b)}{\Gamma(\cb)\Gamma(1-\d)}
-\sum\limits_{k=3}^{m'}\frac{A_{k}}{b_{k}}
{}_{p+1}F_{p}\!\left(\!\left.\!\begin{array}{l}1-\a+b_{k},1-\cb+b_{k},b_{k}\\1-\b_{[k]}+b_{k},1-\d+b_{k},b_{k}+1\end{array}\!\right|(-1)^{\varkappa}\!\right)
\\
-\sum\limits_{k=1}^{n}\frac{B_{k}}{(1-a_{k})}
{}_{p+1}F_{p}\!\left(\!\left.\!\begin{array}{l}1+\b-a_{k},1+\d-a_{k},1-a_{k}\\1+\a_{[k]}-a_{k},1+\cb-a_{k},2-a_{k}\end{array}\!\right|(-1)^{\varkappa}\!\right),
\end{multline*}
where $\varkappa=p-m'-n$ and the numbers $A_k$, $B_k$ are defined in \eqref{eq:constants}
with $b_2=b_1+m$ on the right hand side.  Simple transformations of the above equality lead to 
\begin{multline*}
\sum\limits_{k=0}^{\infty}\frac{[(-1)^{\varkappa}]^k(1-\a+b_{2})_{k}(1-\cb+b_{2})_{k}\Psi}{(1-\b_{[1,2]}+b_{2})_{k}(1-\d+b_{2})_{k}(b_{2}+k)(m+k)!k!}
\\
=-\frac{(-1)^{m\varkappa}(1-\b_{[1,2]}+b_{1})_{m}(1-\d+b_{1})_{m}}{(1-\cb+b_1)_{m}(1-\a+b_{1})_{m}}
\\
\times\!\sum\limits_{k=0}^{m-1}\frac{(-1)^{m-k}(m-k-1)![(-1)^{\varkappa}]^k(1-\a+b_1)_{k}(1-\cb+b_{1})_{k}}{(1-\b_{[1,2]}+b_{1})_{k}(1-\d+b_{1})_{k}(b_1+k)k!}
	\\
+\frac{\pi}{b_{2}m!}\Biggl[\sum\limits_{b\in\b_{[1,2]}}\!\!\!\cot(\pi(b-b_{1}))-\sum\limits_{c\in\cb}\!\!\cot(\pi(c-b_{1}))\Biggr]
\\
\times{}_{p+1}F_{p}\!\left(\!\left.\!\begin{array}{l}1-\a+b_{2},1-\cb+b_{2},b_{2}\\1+m,1-\b_{[1,2]}+b_{2},1-\d+b_{2},b_{2}+1\end{array}\!\right|(-1)^{\varkappa}\!\right)	
\\	
+\frac{(-1)^m\Gamma(\cb-b_{2})\Gamma(1-\d+b_{2})}{\Gamma(\b_{[1,2]}-b_{2})\Gamma(1-\a+b_{2})}\Bigg\{
\frac{\Gamma(1-\a)\Gamma(\b)}{\Gamma(\cb)\Gamma(1-\d)}
-\sum\limits_{k=3}^{m'}
\frac{\Gamma(\b_{[k]}-b_{k})\Gamma(1-\a+b_k)}{b_{k}\Gamma(\cb-b_{k})\Gamma(1-\d+b_{k})}
\\
\times{}_{p+1}F_{p}\!\left(\!\left.\!\begin{array}{l}1-\a+b_{k},1-\cb+b_{k},b_{k}\\1-\b_{[k]}+b_{k},1-\d+b_{k},b_{k}+1\end{array}\!\right|(-1)^{\varkappa}\!\right)
\\
-\!\sum\limits_{k=1}^{n}\!
\frac{\Gamma(a_k-\a_{[k]})\Gamma(1+\b-a_{k})}{(1\!-\!a_{k})\Gamma(a_{k}\!-\!\d)\Gamma(1\!+\!\cb\!-\!a_{k})}
{}_{p+1}F_{p}\!\!\left(\!\!\left.\!\begin{array}{l}1+\b-a_{k},1+\d-a_{k},1-a_{k}\\1+\a_{[k]}-a_{k},1+\cb-a_{k},2-a_{k}\end{array}\!\!\right|\!(-1)^{\varkappa}\!\right)\!\!\Bigg\}.
\end{multline*}
Denote
$$
\aalpha=1-\a+b_1+m,~~\bbeta=1-\cb+b_1+m,~~\ggamma=1-\b_{[1,2]}+b_1+m,~~\ddelta=1-\d+b_1+m
$$
and recall that $b_2=b_1+m$, so that
$$
\a=1-\aalpha+b_2,~~\b=(b_2-m,b_2,1-\ggamma+b_2),~~\cb=1-\bbeta+b_2,~~\d=1-\ddelta+b_2.
$$
Then, the latter identity takes the form (recall that $b_2=b_1+m$):
\begin{multline*}
\frac{1}{b_{2}m!}\sum\limits_{k=0}^{\infty}\frac{(\aalpha)_{k}(\bbeta)_{k}(b_{2})_{k}[(-1)^{\varkappa}]^{k}}{(\ggamma)_{k}(\ddelta)_{k}
(b_{2}+1)_{k}(1+m)_{k}k!}\Psi_1
\\
=\frac{\pi}{m!b_{2}}
\left[\sum\limits_{\beta\in\bbeta}\cot(\pi\beta)-\sum\limits_{\gamma\in\ggamma}\cot(\pi\gamma)\right]
{}_{p+1}F_{p}\!\left(\!\left.\!\begin{array}{l}\aalpha,\bbeta,b_{2}\\1+m,\ggamma,\ddelta,b_{2}+1\end{array}\!\right|(-1)^{\varkappa}\!\right)
\\
+\frac{(-1)^{m}\Gamma(1-\bbeta)\Gamma(\ddelta)}{\Gamma(1-\ggamma)\Gamma(\aalpha)}
\Bigg\{\frac{\Gamma(\aalpha-b_{2})\Gamma(b_2-m)\Gamma(b_2)\Gamma(1-\ggamma+b_{2})}{\Gamma(1-\bbeta+b_2)\Gamma(\ddelta-b_2)}
\\
-\sum\limits_{k=1}^{m'-2}\frac{\Gamma(\gamma_k-m-1)\Gamma(\gamma_k-1)\Gamma(\gamma_{k}-\ggamma_{[k]})\Gamma(1-\gamma_k+\aalpha)}{(1-\gamma_{k}+b_{2})\Gamma(\gamma_k-\bbeta)\Gamma(1-\gamma_k+\ddelta)}
\\
\times{}_{p+1}F_{p}\!\left(\!\left.\!\begin{array}{l}1-\gamma_k+\aalpha,1-\gamma_{k}+\bbeta,1-\gamma_{k}+b_{2}\\2+m-\gamma_k,2-\gamma_k,
1-\gamma_{k}+\ggamma_{[k]},1-\gamma_k+\ddelta,2-\gamma_{k}+b_{2}\end{array}\!\right|(-1)^{\varkappa}\!\right)
\\
-\sum\limits_{k=1}^{n}\frac{\Gamma(\aalpha_{[k]}-\alpha_{k})\Gamma(\alpha_k-m)\Gamma(\alpha_k)\Gamma(1-\ggamma+\alpha_k)}{(\alpha_{k}-b_{2})\Gamma(\ddelta-\alpha_k)\Gamma(1-\bbeta+\alpha_{k})}
\\
\times {}_{p+1}F_{p}\!\left(\!\left.\!\begin{array}{l}\alpha_{k}-m,\alpha_{k},1-\ggamma+\alpha_{k},1-\ddelta+\alpha_{k},\alpha_{k}-b_{2}
\\
1-\aalpha_{[k]}+\alpha_{k},1-\bbeta+\alpha_{k},1+\alpha_{k}-b_{2}\end{array}\!\right|(-1)^{\varkappa}\!\right)\Bigg\}
\\
-\frac{(1-\ggamma)_{m}(1-\ddelta)_{m}}{(1-\aalpha)_{m}(1-\bbeta)_{m}}
\sum\limits_{k=0}^{m-1}\frac{(m-k-1)![(-1)^{\varkappa-1}]^{k+m}(\aalpha-m)_{k}(\bbeta-m)_{k}}{(\ggamma-m)_{k}(\ddelta-m)_{k}(b_2-m+k)k!},
\end{multline*}
where
\begin{multline*}
\Psi_1=\sum\limits_{\gamma\in\ggamma}\!\psi(\gamma+k)+\sum\limits_{\delta\in\ddelta}\!\psi(\delta+k)+\psi(1+m+k)+\psi(1+k)
\\
-\!\sum\limits_{\alpha\in\aalpha}\!\psi(\alpha+k)-\!\sum\limits_{\beta\in\bbeta}\!\psi(\beta+k)+\frac{1}{b_2+k}.
\end{multline*}
In view of \eqref{eq:m-mtrick} we have
\begin{multline*}
\frac{(1-\ggamma)_{m}(1-\ddelta)_{m}}{(1-\aalpha)_{m}(1-\bbeta)_{m}}
\sum\limits_{k=0}^{m-1}\frac{(m-k-1)![(-1)^{\varkappa-1}]^{k+m}(\aalpha-m)_{k}(\bbeta-m)_{k}}{(\ggamma-m)_{k}(\ddelta-m)_{k}(b_2-m+k)k!}
\\
=\frac{(-1)^{\varkappa-1}(1-\ggamma)_{1}(1-\ddelta)_{1}}{(m-1)!(1-\aalpha)_{1}(1-\bbeta)_{1}}
\sum\limits_{j=0}^{m-1}\frac{(1-m)_j(2-\ggamma)_{j}(2-\ddelta)_{j}[(1)_j]^2(-1)^{\varkappa{j}}}{(2-\bbeta)_{j}(2-\aalpha)_{j}(b_2-1-j)j!}
\\
=\!\frac{(-1)^{\varkappa-1}(1\!-\!\ggamma)_{1}(1\!-\!\ddelta)_{1}}{(m\!-\!1)!(1\!-\!\aalpha)_{1}(1\!-\!\bbeta)_{1}(b_2\!-\!1)}
{}_{p+2}F_{p+1}\!\left(\!\!\left.\!\begin{array}{l}1-m,2-\ggamma,2-\ddelta,1,1,1-b_2
\\
2-\aalpha,2-\bbeta,2-b_2\end{array}\!\!\right|\!(-1)^{\varkappa}\!\right).
\end{multline*}

Next, we put $\aalpha=(1+m,b_{2}+1,\aalpha')$ and substitute into the large identity above:
\begin{multline*}
\sum\limits_{k=0}^{\infty}\frac{(\aalpha')_{k}(\bbeta)_{k}(b_{2})_{k}[(-1)^{\varkappa}]^{k}}{(\ggamma)_{k}(\ddelta)_{k}
k!}
\\
\times\!\bigg\{\!\sum\limits_{\gamma\in\ggamma}\psi(\gamma+k)+\sum\limits_{\delta\in\ddelta}\psi(\delta+k)+\psi(1+k)-\sum\limits_{\alpha\in\aalpha'}\psi(\alpha+k)-\sum\limits_{\beta\in\bbeta}\psi(\beta+k)-\psi(b_{2}+k)\!\bigg\}
\\
=\pi
\left[\sum\limits_{\beta\in\bbeta}\cot(\pi\bbeta)-\sum\limits_{\gamma\in\ggamma}\cot(\pi\ggamma)\right]
{}_{p-1}F_{p-2}\!\left(\!\left.\!\begin{array}{l}\aalpha',\bbeta,b_{2}\\\ggamma,\ddelta\end{array}\!\right|(-1)^{\varkappa}\!\right)
\\
+\frac{(-1)^{m}\Gamma(1-\bbeta)\Gamma(\ddelta)}{\Gamma(1-\ggamma)\Gamma(\aalpha')\Gamma(b_{2})}
\Bigg\{\frac{\pi\Gamma(\aalpha'-b_{2})\Gamma(b_2)\Gamma(1-\ggamma+b_{2})}{(-1)^{m}\sin(\pi{b_2})\Gamma(1-\bbeta+b_2)\Gamma(\ddelta-b_2)}
\\
\times\!{}_{p-1}F_{p-2}\!\left(\!\left.\!\begin{array}{l}b_{2},1-\ggamma+b_{2},1-\ddelta+b_{2}
\\
1-\aalpha'+b_{2},1-\bbeta+b_{2}\end{array}\!\right|(-1)^{\varkappa}\!\right)
\\
-\sum\limits_{k=1}^{m'-2}\frac{\pi\Gamma(\gamma_k-1)\Gamma(\gamma_{k}-\ggamma_{[k]})\Gamma(1-\gamma_k+\aalpha')\Gamma(1-\gamma_k+b_{2})}{(-1)^m\sin(\pi(\gamma_k-1))\Gamma(\gamma_k-\bbeta)\Gamma(1-\gamma_k+\ddelta)}
\\
\times {}_{p-1}F_{p-2}\!\left(\!\left.\!\begin{array}{l}1-\gamma_k+\aalpha',1-\gamma_{k}+\bbeta,1-\gamma_{k}+b_{2}\\2-\gamma_k,
1-\gamma_{k}+\ggamma_{[k]},1-\gamma_k+\ddelta\end{array}\!\right|(-1)^{\varkappa}\!\right)
\\
+\sum\limits_{k=1}^{n-2}\frac{\pi\Gamma(b_{2}-\alpha'_{k})\Gamma(\aalpha'_{[k]}-\alpha'_{k})\Gamma(\alpha'_k)\Gamma(1-\ggamma+\alpha'_k)}{(-1)^m\sin(\pi\alpha'_k)\Gamma(\ddelta-\alpha'_k)\Gamma(1-\bbeta+\alpha'_{k})}
\\
\times\! {}_{p-1}F_{p-2}\!\left(\!\left.\!\begin{array}{l}\alpha'_{k},1-\ggamma+\alpha'_{k},1-\ddelta+\alpha'_{k}
\\
1-b_{2}+\alpha'_{k},1-\aalpha'_{[k]}+\alpha_{k},1-\bbeta+\alpha'_{k}\end{array}\!\right|(-1)^{\varkappa}\!\right)\Bigg\}
\\
-\frac{(-1)^{\varkappa}(1-\ggamma)_{1}(1-\ddelta)_{1}}{(1-b_{2})(1-\aalpha')_{1}(1-\bbeta)_{1}}
{}_{p}F_{p-1}\!\left(\!\left.\!\begin{array}{l}1,1,2-\ggamma,2-\ddelta\\
2-b_{2},2-\aalpha',2-\bbeta\end{array}\!\right|(-1)^{\varkappa}\!\right),
\end{multline*}
where we applied the easy verifiable identities
\begin{multline*}
\frac{(-1)^{m}\Gamma(1-\bbeta)\Gamma(\ddelta)}{\Gamma(1-\ggamma)\Gamma(\aalpha')\Gamma(b_{2})}
\frac{\Gamma(\aalpha'-m-1)\Gamma(b_2-m)m!\Gamma(2-\ggamma+m)}{(m+1-b_{2})\Gamma(\ddelta-m-1)\Gamma(2-\bbeta+m)}
\\
\times{}_{p}F_{p-1}\!\left(\!\left.\!\begin{array}{l}1,m+1,2-\ggamma+m,2-\ddelta+m
\\
2-\aalpha'+m,2-\bbeta+m,2-b_{2}+m\end{array}\!\right|(-1)^{\varkappa}\!\right)
\\
+\frac{(-1)^{m}(1-\ggamma)_{m}(1-\ddelta)_{m}}{(1-b_{2})_{m}(1-\aalpha')_{m}(1-\bbeta)_{m}}
\sum\limits_{k=0}^{m-1}\frac{(m-k-1)!(1)_{k}(b_{2}-m)_{k}(\aalpha'-m)_{k}(\bbeta-m)_{k}}{[(-1)^{\varkappa-1}]^{k+m}(\ggamma-m)_{k}(\ddelta-m)_{k}k!}
\\
=\frac{(-1)^{\varkappa}(1-\ggamma)_{1}(1-\ddelta)_{1}}{(1-b_{2})(1-\aalpha')_{1}(1-\bbeta)_{1}}
{}_{p}F_{p-1}\!\left(\!\left.\!\begin{array}{l}1,1,2-\ggamma,2-\ddelta\\
2-b_{2},2-\aalpha',2-\bbeta\end{array}\!\right|(-1)^{\varkappa}\!\right)
\end{multline*}
and
\begin{multline*}
\frac{\pi\Gamma(\aalpha'-b_{2})\Gamma(b_2)\Gamma(1-\ggamma+b_{2})}{(-1)^{m}\sin(\pi{b_2})\Gamma(1\!-\!\bbeta\!+\!b_2)\Gamma(\ddelta\!-\!b_2)}
+\frac{\pi\Gamma(\aalpha'-b_{2}-1)\Gamma(b_{2}+1)\Gamma(2-\ggamma+b_{2})}{(-1)^{m}\sin(\pi{b_{2}})\Gamma(\ddelta-b_{2}-1)\Gamma(2-\bbeta+b_{2})}
\\
\times{}_{p}F_{p-1}\!\left(\!\left.\!\begin{array}{l}b_{2}+1,2-\ggamma+b_{2},2-\ddelta+b_{2},1
\\
2-\aalpha'+b_{2},2-\bbeta+b_{2},2\end{array}\!\right|(-1)^{\varkappa}\!\right)
\\
=\!\frac{\pi\Gamma(\aalpha'\!-\!b_{2})\Gamma(b_2)\Gamma(1-\ggamma+b_{2})}{(-1)^{m}\sin(\pi{b_2})\Gamma(1\!-\!\bbeta\!+\!b_2)\Gamma(\ddelta\!-\!b_2)}
{}_{p-1}F_{p-2}\!\left(\!\left.\!\!\begin{array}{l}b_{2},1-\ggamma+b_{2},1-\ddelta+b_{2}
\\
1-\aalpha'+b_{2},1-\bbeta+b_{2}\end{array}\!\right|(-1)^{\varkappa}\!\right).
\end{multline*}

Hence, writing $\aalpha=(b_2,\aalpha')\in\C^{n-1}$ and recalling that $\bbeta\in\C^{p-n}$, $\ggamma\in\C^{m'-2}$, $\ddelta\in\C^{p-m'}$, $\varkappa=p-n-m'$ and $m'+n\ge{p}$, we get
\begin{multline*}
\sum\limits_{k=0}^{\infty}\frac{(\aalpha)_{k}(\bbeta)_{k}[(-1)^{\varkappa}]^{k}}{(\ggamma)_{k}(\ddelta)_{k}
k!}
\bigg\{\!\sum\limits_{\gamma\in\ggamma}\psi(\gamma+k)+\sum\limits_{\delta\in\ddelta}\psi(\delta+k)+\psi(1+k)\!-\!\sum\limits_{\alpha\in\aalpha}\psi(\alpha+k)
\\
-\sum\limits_{\beta\in\bbeta}\psi(\beta+k)\!\bigg\}
=\pi\left[\sum\limits_{\beta\in\bbeta}\cot(\pi\beta)-\sum\limits_{\gamma\in\ggamma}\cot(\pi\gamma)\right]
{}_{p-1}F_{p-2}\!\left(\!\left.\!\begin{array}{l}\aalpha,\bbeta\\\ggamma,\ddelta\end{array}\!\right|(-1)^{\varkappa}\!\right)
\\
-\frac{(-1)^{\varkappa}(1-\ggamma)_{1}(1-\ddelta)_{1}}{(1-\aalpha)_{1}(1-\bbeta)_{1}}
{}_{p}F_{p-1}\!\left(\!\left.\!\begin{array}{l}1,1,2-\ggamma,2-\ddelta\\
2-\aalpha,2-\bbeta\end{array}\!\right|(-1)^{\varkappa}\!\right)
\\
+\frac{\Gamma(1-\bbeta)\Gamma(\ddelta)}{\Gamma(1-\ggamma)\Gamma(\aalpha)}\sum\limits_{k=1}^{n-1}\frac{\pi\Gamma(\aalpha_{[k]}-\alpha_{k})\Gamma(\alpha_k)\Gamma(1-\ggamma+\alpha_k)}{\sin(\pi\alpha_k)\Gamma(\ddelta-\alpha_k)\Gamma(1-\bbeta+\alpha_{k})}
\\
\times\!{}_{p-1}F_{p-2}\!\left(\!\left.\!\begin{array}{l}\alpha_{k},1-\ggamma+\alpha_{k},1-\ddelta+\alpha_{k}
\\
1-\aalpha_{[k]}+\alpha_{k},1-\bbeta+\alpha_{k}\end{array}\!\right|(-1)^{\varkappa}\!\right)
\\
+\frac{\Gamma(1-\bbeta)\Gamma(\ddelta)}{\Gamma(1-\ggamma)\Gamma(\aalpha)}\sum\limits_{k=1}^{m'-2}\frac{\pi\Gamma(\gamma_k-1)\Gamma(\gamma_{k}-\ggamma_{[k]})\Gamma(1-\gamma_k+\aalpha)}{\sin(\pi\gamma_k)\Gamma(\gamma_k-\bbeta)\Gamma(1-\gamma_k+\ddelta)}
\\
\times\!{}_{p-1}F_{p-2}\!\left(\!\left.\!\begin{array}{l}1-\gamma_k+\aalpha,1-\gamma_{k}+\bbeta\\2-\gamma_k,
1-\gamma_{k}+\ggamma_{[k]},1-\gamma_k+\ddelta\end{array}\!\right|(-1)^{\varkappa}\!\right).
\end{multline*}
Introducing  the notation
\begin{align*}
&\bbeta=\a_1\in\C^{n},~~ \aalpha=\a_2\in\C^{p-n},~~ \a=(\a_1,\a_2)\in\C^{p}, 
\\
&\b_1=\ggamma\in\C^{m},~~ \b_2=\ddelta\in\C^{p-1-m},~~ \b=(\b_1,\b_2)\in\C^{p-1},    
\end{align*}
we finally arrive at \eqref{eq:alternatingFinal}  under the conditions $0\le{n}\le{p}$, $0\le{m}\le{p-1}$ and $\kappa:=m-n\ge-2$. 
$\hfill\square$

\medskip

A particular case of this identity for $n=p$, $m=p-1$ takes the form:
\begin{multline}\label{eq:alternatingFinalMax}
\sum\limits_{k=0}^{\infty}\frac{(\a)_{k}(-1)^{k}}{\pi(\b)_{k}
k!}
\bigg[\sum_{b\in\b}\psi(b+k)+\psi(1+k)-\sum_{a\in\a}\psi(a+k)\bigg]\!=\!
\\
\left[\sum_{a\in\a}\cot(\pi{a})\!-\!\sum_{b\in\b}\cot(\pi{b})\right]
{}_{p}F_{p-1}\!\left(\!\left.\!\begin{array}{l}\a\\\b\end{array}\!\right|-1\!\right)
+\frac{(1-\b)_{1}}{\pi(1-\a)_{1}}
{}_{p+1}F_{p}\!\left(\!\left.\!\begin{array}{l}1,1,2-\b\\
2-\a\end{array}\!\right|-1\!\right)
\\
+\frac{\Gamma(1-\a)}{\Gamma(1-\b)}\sum\limits_{b_k\in\b}
\frac{\Gamma(b_k-1)\Gamma(b_{k}-\b_{[k]})}{\sin({\pi}b_k)\Gamma(b_k-\a)}
{}_{p}F_{p-1}\!\left(\!\left.\!\begin{array}{l}1-b_k+\a\\2-b_k,
1-b_{k}+\b_{[k]}\end{array}\!\right|-1\!\right).
\end{multline}

Combining this formula with \eqref{eq:final} we obtain:
\begin{multline*}
\sum\limits_{k=0}^{\infty}\frac{(\a)_{k}(\pm1)^{k}}{\pi(\b)_{k}
k!}
\bigg[\sum_{b\in\b}\psi(b+k)+\psi(1+k)-\sum_{a\in\a}\psi(a+k)\bigg]
\\
\!=\!\!\left[\sum_{a\in\a}\cot(\pi{a})\!-\!\sum_{b\in\b}\cot(\pi{b})\right]
\!{}_{p}F_{p-1}\!\!\left(\!\!\left.\!\begin{array}{l}\a\\\b\end{array}\!\right|\pm1\!\right)
\mp\frac{(1-\b)_{1}}{\pi(1-\a)_{1}}
{}_{p+1}F_{p}\!\!\left(\!\!\left.\!\begin{array}{l}1,1,2-\b\\
2-\a\end{array}\!\right|\pm1\!\right)
\\
+\frac{\Gamma(1-\a)}{\Gamma(1-\b)}\sum\limits_{b_k\in\b}
\!\frac{[\cos({\pi}b_k)]^{\genfrac\{\}{0pt}{}{1}{0}}\Gamma(b_k-1)\Gamma(b_{k}-\b_{[k]})}{\sin({\pi}b_k)\Gamma(b_k-\a)}
{}_{p}F_{p-1}\!\!\left(\!\!\left.\!\begin{array}{l}1-b_k+\a\\2-b_k,
1-b_{k}+\b_{[k]}\end{array}\!\right|\pm1\!\right)
\end{multline*}
which is precisely formula \eqref{eq:FinalMax} given in the Introduction. Formula  \eqref{eq:FinalMin}, on the other hand, is obtained by choosing $m=n=0$ 
in \eqref{eq:alternatingFinal}:
\begin{multline*}
\sum\limits_{k=0}^{\infty}\frac{(\a)_{k}}{\pi(\b)_{k}
k!}
\bigg[\sum_{b\in\b}\psi(b+k)+\psi(1+k)
-\sum_{a\in\a}\psi(a+k)\bigg]
\\
\!=\!\frac{\Gamma(\b)}{\Gamma(\a)}
\sum\limits_{a_k\in\a}\frac{\Gamma(a_k)\Gamma(\a_{[k]}-a_{k})}{\sin(\pi a_k)\Gamma(\b-a_k)}
{}_{p}F_{p-1}\!\left(\!\!\begin{array}{l}a_{k},1-\b+a_{k}
\\
1-\a_{[k]}+a_{k}\end{array}\!\!\right)
-\frac{(1-\b)_{1}}{\pi(1-\a)_{1}}
{}_{p+1}F_{p}\!\left(\!\!\begin{array}{l}1,1,2-\b\\
2-\a\end{array}\!\!\right).
\end{multline*}
Combining it with the "+" case of \eqref{eq:FinalMax} we also get:
\begin{multline}\label{eq:FinalNoPsi}
\left[\sum_{a\in\a}\cot(\pi{a})\!-\!\sum_{b\in\b}\cot(\pi{b})\right]
{}_{p}F_{p-1}\!\left(\!\!\begin{array}{l}\a\\\b\end{array}\!\right)
\\
=\frac{\Gamma(\b)}{\Gamma(\a)}
\sum\limits_{k=1}^{p}\frac{\Gamma(a_k)\Gamma(\a_{[k]}-a_{k})}{\sin(\pi a_k)\Gamma(\b-a_k)}
{}_{p}F_{p-1}\!\left(\!\!\begin{array}{l}a_{k},1-\b+a_{k}
\\
1-\a_{[k]}+a_{k}\end{array}\!\!\right)
\\
-\frac{\Gamma(1-\a)}{\Gamma(1-\b)}\sum\limits_{k=1}^{p-1}
\cot({\pi}b_k)\frac{\Gamma(b_k-1)\Gamma(b_{k}-\b_{[k]})}{\Gamma(b_k-\a)}
{}_{p}F_{p-1}\!\left(\!\!\begin{array}{l}1-b_k+\a\\2-b_k,
1-b_{k}+\b_{[k]}\end{array}\!\right)
\end{multline}
- another multi-term hypergeometric identity which appears to be new and combines the feature of \cite[(5.1)]{CKP} and \cite[(5.3)]{CKP}.

\end{document}